\documentclass[12pt,reqno,?]{amsart2000}    
\usepackage{dcolumn}
\usepackage{epic}
\usepackage{eepicemu}
\usepackage{float}
\allowdisplaybreaks[4]                 
\theoremstyle{plain}
\newtheorem*{lemmaoneone}{Lemma I$\boldsymbol{{}^{(1)}}$}
\newtheorem*{lemmaonetwo}{Lemma I$\boldsymbol{{}^{(2)}}$}

\newtheorem*{lemmatwotwo}{Lemma II$\boldsymbol{{}^{(2)}}$}
\newtheorem*{lemmaone}{Lemma I}
\newtheorem*{lemmatwo}{Lemma II}
\newtheorem*{untitled}{}

\theoremstyle{definition}
\newtheorem*{example}{Example}
\newtheorem*{normalization}{Normalization}
\newtheorem*{notation}{Notation}

\newcommand{\tn}{\boldsymbol{\mathsf{T}}}
\newcommand{\sn}{\boldsymbol{\mathsf{S}}}  
\newcommand{\pn}{\boldsymbol{\mathsf{P}}}  
\newcommand{\mn}{\boldsymbol{\mathsf{M}}}  
\newcommand{\cn}{\boldsymbol{\mathsf{C}}}  
\newcommand{\fn}{\boldsymbol{\mathsf{F}}}  
\newcommand{\gfn}{\fn^{\displaystyle{*}}}  
\DeclareMathOperator{\type}{type}
\newcommand{\tstrut}{\rule[-1.25ex]{0ex}{3.75ex}}
\newcommand{\h}{$\phantom{0}$}
\newcommand{\K}{\h\h}
\newcommand{\timeno}{2001-0924}
\newcommand{\pageno}{29}
\newcommand{\header}{\textsc{Structure of Binary Sequences 
                     \copyright\ J.~Tharrats 1982--2001 Printed \timeno} 
                     \texttt{math.CO/0109182}}
\begin{document}
%
%
\title[\protect\header]{Structure of Binary Sequences}
\date{February, 1982}
\thanks{\emph{Report No.}: UPR-PP-82/1; \copyright\ J.~Tharrats, 1982--2001; 
        \emph{Printed}: \timeno; \pageno\ pp}
\thanks{\emph{arXiv e-print}\texttt{ math.CO/0109182 }}
%
%
\author[\protect\header]{J.~Tharrats}
\address{Department of Physics, University of Puerto Rico at R\'{\i}o Piedras, 
         PO BOX 23343, SAN JUAN PR 00931-3343} 
\email{tharrats@rrpac.upr.clu.edu}
%
%
\begin{abstract}
The distribution of a given sequence in the set of all sequences with 
$n$ ones and $m = M - n$ zeros are found by relating the problem to the 
partitions of a natural number in $m$ natural summands, taking into
account the order. The formulas obtained have many applications both
in Physics and Mathematics. Examples discussed in the present paper
are: non Markovian chains, partition functions of binary alloys and
Ising magnets, generalized Kaplansky lemma, generalized Fibonacci
numbers and a general expansion of 
$\sum_{h=0}^{m} h^{r} {\binom{m}{h}}^{2}$ 
in terms of the Stirling numbers of second kind.
\end{abstract}
\maketitle
\section{Introduction}\label{s-I} 

Some combinatorial problems can be solved immediately by refering them to the 
structure of the $\binom{m+n}{n}$ sequences of $m$ zeros and $n$ ones of length 
$N=m+n$, i.e.: knowing the ocurrence distribution of a given ordered string of 
$r$ $0$'s and $s$ $1$'s ($r+s < N$). Take, for example the problem of finding 
the number of subsets of cardinality $n$ of the set $X = \{1,2,\ldots,N\}$ that 
contain exactly $h$ times a number $r$ of consecutive naturals; this is done 
knowing how many times the string $(1,1,\stackrel{r}{\ldots},1)$ has $h$ 
recurrences in those sequences. For a non fixed cardinality one makes the sum 
over $n$ and obtains how many elements of the power set have the said property. 
In this case, as one can see in Riordan \cite[problem 1(b), p.\ 14]{Rj1958}, 
for $h=0$, $r=2$ one obtains the Fibonacci numbers, though Carlitz 
\cite{Cl1974} had generalized to $h=0$, $r > 2$; (in 1976 Carlitz \cite{Cl1976} 
precises the problem mixing strings of $(00)$ and $(11)$).

Thus, the structure of digit-sequences can have its own place in
Combinatorics, introducing problems concerning the sequences
themselves; for instance, in the $\binom{N}{n}$ sequences, the
number of sequences that have each pair of ones separated by at least
$k$ zeros. Although the solution is given by Kaplansky lemma
\cite{Ki1953}, the problem becomes trivial by applying the formula
which gives the distribution of the $(0,0,\stackrel{k}{\ldots},0;1)$
strings, choosing the result which gives exactly $n$ strings. More
generally, we will also find the case in which only $s$ ones will be
separated by at least $k$ zeros.

In this theory two points of view can be taken: One may decide to
bound the sequences with first and last digit or to consider the first
digit following the last in cyclic fashion. Since our main
interests are physical applications, only the second point of view
will be followed. This view, called here ``Ising process'' gives in
the particular case above the corrected Fibonacci numbers. It is also
more symmetrical, which explains why our method is so simple and can
be developed free of generating functions, in contrast with the method
of bounded sequences in which the particular cases solved by Carlitz
give results strongly dependent on the first digit.

This paper was motivated by the solution of non-markovian random 
walks.\footnote{The paper ``Analytical Solution of non-Markovian Random Walks'' 
in which we give some applications to crystal lattices, is to appear.} As is 
well-known, the solution of the one-dimensional random walk, with 
probabilities $p$ (moving to the right) and $q = 1 - p$ (to the left), 
is given by Bernouilli's distribution:
\begin{equation}\label{e-1}
    P(k) 
    = \binom{N}{\frac{1}{2}(N+k)} p^{\frac{N+k}{2}} q^{\frac{N-k}{2}}.
\end{equation}
But if transition probabilities are dependent on past states, i.e., if
the particle has memory, then the distribution changes radically. In
the simple case of dependence on the last state only, the probability
is $\alpha p$ or $\alpha q$ for particle moving backwards and $\beta p$ 
or $\beta q$ moving forwards ($\alpha + \beta = 1$). Let
$\tn^{mn}_{\tau}$ be the number of paths of length $N$ with
$\tau$ changes in direction of the particle, or ``jumps'' in $N$ steps
(number of ocurrences of $(01)$ or $(10)$). Then paths like $ppqqp$ and
$pqpqp$ have different probability because they are affected by the 
factors $\alpha^{2}\beta^{3}$ and $\alpha^{4}\beta$ respectively, and 
the binomial coefficient in \eqref{e-1} splits in the form
\[
    \binom{N}{\frac{1}{2} (N+k)} \longrightarrow 
    \sum_{\tau=0}^{N-k} \tn^{\frac{N+k}{2},\frac{N-k}{2}}_{\tau}
                        \alpha^{\tau} \beta^{N-\tau}    .
\]
Thus the formula \eqref{e-1} changes depending on the constant $\alpha$ 
which measures the tendency to go forwards ($\alpha < \frac{1}{2}$) or
backwards ($\alpha > \frac{1}{2}$).  We see that the problem is solved
by finding the $\tn$-numbers.

We also will see how $\tn$-numbers give the natural way to find some
Thermodynamic Partition Functions, such as in Ising magnets, in
contrast with naive ways given by some physicists.

Historically, it seems Andr\'e \cite{Ad1910} was the first to consider
such sequences (see Comptet \cite{Cl1970}) in order to solve
Bertrand's scrutiny problem. He considered minimal paths in a square
lattice where $y$-jumps stand for $0$'s and $x$-jumps for $1$'s. This
give a simple way to solve Bertrand's problem counting the number of
paths which do not touch the $y=x$ line. The solution is not
invariant by a circular permutation of the sequences and falls in the
case of unbounded sequences, for the path $(0010001)$ is allowable (in
each step the eventual winner is winning) but the Ising equivalent
path $(0100010)$ is not. Processes like this one are also non-markovian
because allowed decisions depend of all anterior states, but can be
solved by elementary methods; the scrutiny problem can be solved as in
\cite{Cl1970} or by reducing it to well known random walks avoiding
the origin (projecting paths on $y=-x$).

Along the way in this paper many combinatorial formulas arise, some
very well known but others related to Stirling numbers are new, as far
as we know; but what really matters is the natural way in which they 
arise.

\begin{notation}
We prefer to preserve the (old fashioned?) symbols for: 
$\cn^{m}_{n} = \binom{n}{m}$; $\sn^{m}_{n} =$ Stirling numbers of the 
second kind; $\pn^{m}_{n} =$ number of partitions of a natural number 
$n$ in $m$ natural summands disregarding the order; $\mn^{m}_{n}$ for
number of partitions regarding the order ($\mn$ stands for De Moivre 
who discovered it to be equal to $\binom{n-1}{m-1}$).
\end{notation}

As $\pn$-partitions can be associated with Young-Sylvester's tableaux,
$\mn$-partitions too, taking into account the order.  For example, the
partition of 6 into 3 parts gives the $\pn^{3}_{6} = 3$ tableaux:
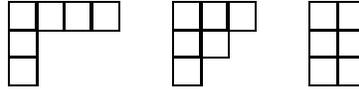
\begin{center}
\setlength{\unitlength}{2ex}
\begin{picture}(13,3)
%
\multiput( 0,2)(1,0){4}{\framebox(1,1){}}
\multiput( 0,1)(1,0){1}{\framebox(1,1){}}
\multiput( 0,0)(1,0){1}{\framebox(1,1){}}
\multiput( 6,2)(1,0){3}{\framebox(1,1){}}
\multiput( 6,1)(1,0){2}{\framebox(1,1){}}
\multiput( 6,0)(1,0){1}{\framebox(1,1){}}
\multiput(11,2)(1,0){2}{\framebox(1,1){}}
\multiput(11,1)(1,0){2}{\framebox(1,1){}}
\multiput(11,0)(1,0){2}{\framebox(1,1){}}
\end{picture}
\end{center}
For $\mn$-partitions, we have the following
$\mn^{3}_{6} = \binom{5}{2} = 
\frac{3!}{2!\,1!} + \frac{3!}{1!\,1!\,1!} + \frac{3!}{3!} = 10$ 
tableaux:
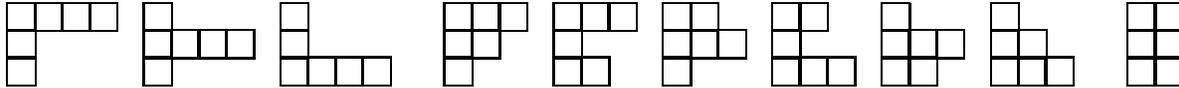
\begin{center}
\setlength{\unitlength}{2ex}
\begin{picture}(43,3)
%
\multiput( 0,2)(1,0){4}{\framebox(1,1){}}
\multiput( 0,1)(1,0){1}{\framebox(1,1){}}
\multiput( 0,0)(1,0){1}{\framebox(1,1){}}
\multiput( 5,2)(1,0){1}{\framebox(1,1){}}
\multiput( 5,1)(1,0){4}{\framebox(1,1){}}
\multiput( 5,0)(1,0){1}{\framebox(1,1){}}
\multiput(10,2)(1,0){1}{\framebox(1,1){}}
\multiput(10,1)(1,0){1}{\framebox(1,1){}}
\multiput(10,0)(1,0){4}{\framebox(1,1){}}
\multiput(16,2)(1,0){3}{\framebox(1,1){}}
\multiput(16,1)(1,0){2}{\framebox(1,1){}}
\multiput(16,0)(1,0){1}{\framebox(1,1){}}
\multiput(20,2)(1,0){3}{\framebox(1,1){}}
\multiput(20,1)(1,0){1}{\framebox(1,1){}}
\multiput(20,0)(1,0){2}{\framebox(1,1){}}
\multiput(24,2)(1,0){2}{\framebox(1,1){}}
\multiput(24,1)(1,0){3}{\framebox(1,1){}}
\multiput(24,0)(1,0){1}{\framebox(1,1){}}
\multiput(28,2)(1,0){2}{\framebox(1,1){}}
\multiput(28,1)(1,0){1}{\framebox(1,1){}}
\multiput(28,0)(1,0){3}{\framebox(1,1){}}
\multiput(32,2)(1,0){1}{\framebox(1,1){}}
\multiput(32,1)(1,0){3}{\framebox(1,1){}}
\multiput(32,0)(1,0){2}{\framebox(1,1){}}
\multiput(36,2)(1,0){1}{\framebox(1,1){}}
\multiput(36,1)(1,0){2}{\framebox(1,1){}}
\multiput(36,0)(1,0){3}{\framebox(1,1){}}
\multiput(41,2)(1,0){2}{\framebox(1,1){}}
\multiput(41,1)(1,0){2}{\framebox(1,1){}}
\multiput(41,0)(1,0){2}{\framebox(1,1){}}
\end{picture}
\end{center}
which we call Young-De Moivre tableaux.

In the unbounded poset of integers, we call $m \wedge n$ the
$\min{(m,n)}$ and $m \vee n$ the $\max{(m,n)}$; (no confusion arise 
with connectors of formal systems).

We will refer always the known combinatorial formulas to the
monumental ``Ars'' by Knuth \cite{Kde1968}.\footnote{Knuth uses the
Christoffel symbols of first and second kind for Stirling numbers of
first and second kind respectively.}
\section{$\tn$-Numbers}\label{s-TN}

Let us consider the set of all sequences with $n$ ones and $m = N - n$
zeros; an element of this set is, for example (for $N = 24$, $m = 13$):
\begin{equation}\label{e-2}
    000001110011010001100111 . 
\end{equation}
To every sequence we associate a number $\tau$, which is the number of
jumps from $0$ to $1$, or from $1$ to $0$, including extremes, so in 
the case of the sequence \eqref{e-2} we have $\tau = 10$.  Obviously, 
the number $\tau$ is even and it is invariant under any circular 
permutation of the sequence.

The problem here is to find the distribution of $\tau$'s, i.e.: How 
many $\tau$ equal $2$, how many $\tau$ equal $4$, etc.\ in these
$\binom{N}{m}$ sequences.

The following table shows the case $N = 7$ and $m = 3$ with the
corresponding values of $\tau$:
\begin{table}[H]\caption{}\label{t-1}%
\begin{center}%
\begin{tabular}{ | c c c c c c c c c c | }
    \hline
       & $\tau$ &    & $\tau$ &    & $\tau$ &    & $\tau$ &    & $\tau$ \\
    0000111 & 2 & 0011001 & 4 & 0101010 & 6 & 1000101 & 4 & 1010100 & 6 \\
    0001011 & 4 & 0011010 & 4 & 0101100 & 4 & 1000110 & 4 & 1011000 & 4 \\
    0001101 & 4 & 0011100 & 2 & 0110001 & 4 & 1001001 & 4 & 1100001 & 2 \\
    0001110 & 2 & 0100011 & 4 & 0110010 & 4 & 1001010 & 6 & 1100010 & 4 \\
    0010011 & 4 & 0100101 & 6 & 0110100 & 4 & 1001100 & 4 & 1100100 & 4 \\
    0010101 & 6 & 0100110 & 4 & 0111000 & 2 & 1010001 & 4 & 1101000 & 4 \\
    0010110 & 4 & 0101001 & 6 & 1000011 & 2 & 1010010 & 6 & 1110000 & 2 \\
    \hline
\end{tabular}
\end{center}
\end{table}

The table gives $7$, $21$ and $7$ sequences ($7+21+7 = \binom{7}{3}$) 
for the corresponding $\tau$-values $2$, $4$ and $6$ respectively. 
This example illustrates the big difficulty in getting these numbers by 
``brute force'' for large values of $N$, (large means here $N > 10$).

Looking at \eqref{e-2} we see that the number $\tau$ is related by a 
particular partition of $m$, interconnected by another partition of 
$n$, both having the same height, i.e., the number $\tau$ is exactly 
twice the height of these two partitions:
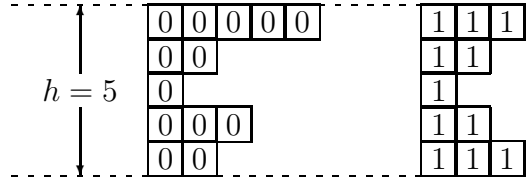
\begin{center}
\setlength{\unitlength}{2.5ex}
\begin{picture}(15,5)(-4,0)
%
\multiput( 0,4)(1,0){5}{\framebox(1,1){0}}
\multiput( 0,3)(1,0){2}{\framebox(1,1){0}}
\multiput( 0,2)(1,0){1}{\framebox(1,1){0}}
\multiput( 0,1)(1,0){3}{\framebox(1,1){0}}
\multiput( 0,0)(1,0){2}{\framebox(1,1){0}}
\multiput( 8,4)(1,0){3}{\framebox(1,1){1}}
\multiput( 8,3)(1,0){2}{\framebox(1,1){1}}
\multiput( 8,2)(1,0){1}{\framebox(1,1){1}}
\multiput( 8,1)(1,0){2}{\framebox(1,1){1}}
\multiput( 8,0)(1,0){3}{\framebox(1,1){1}}
\dashline{0.2}(-4,5)(0,5)    \dashline{0.2}(5,5)(8,5)
\dashline{0.2}(-4,0)(0,0)    \dashline{0.2}(2,0)(8,0)
\put(-2,3.0){\vector(0,1){2}}  
\put(-2,2.5){\makebox(0,0){$h = 5$}}
\put(-2,2.0){\vector(0,-1){2}}
\end{picture}
\end{center}
then, in this case $\tau = 2h = 10$.

We define a ``loop'' associated to two partitions as the set of $2h$
strings obtained by matching these partitions, the first string
regarded as following the last one, In the example we have:
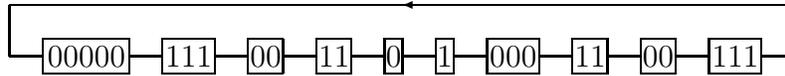
\begin{center}
\setlength{\unitlength}{2.5ex}
\begin{picture}(23,2)
%
\put( 0.0,0.5){\line(1,0){1}} \put( 1.0,0){\framebox(2.5,1){00000}}
\put( 3.5,0.5){\line(1,0){1}} \put( 4.5,0){\framebox(1.5,1){111}}
\put( 6.0,0.5){\line(1,0){1}} \put( 7.0,0){\framebox(1.0,1){00}}
\put( 8.0,0.5){\line(1,0){1}} \put( 9.0,0){\framebox(1.0,1){11}}
\put(10.0,0.5){\line(1,0){1}} \put(11.0,0){\framebox(0.5,1){0}}
\put(11.5,0.5){\line(1,0){1}} \put(12.5,0){\framebox(0.5,1){1}}
\put(13.0,0.5){\line(1,0){1}} \put(14.0,0){\framebox(1.5,1){000}}
\put(15.5,0.5){\line(1,0){1}} \put(16.5,0){\framebox(1.0,1){11}}
\put(17.5,0.5){\line(1,0){1}} \put(18.5,0){\framebox(1.0,1){00}}
\put(19.5,0.5){\line(1,0){1}} \put(20.5,0){\framebox(1.5,1){111}}
\put(22.0,0.5){\line(1,0){1}}
\put( 0.0,0.5){\line(0,1){1.5}} \put( 0.0,2.0){\line( 1,0){11.5}}
\put(23.0,0.5){\line(0,1){1.5}} \put(23.0,2.0){\vector(-1,0){11.5}}
\end{picture}
\end{center}
For a fixed $\tau$, the set of all loops are obtained by matching the
partitions of $m$ with the partitions of $n$, both of height $h =
\frac{1}{2} \tau$. Fortunately this is done not by ordinary partitions
but by \emph{partitions regarding the order}, giving the numbers:
\[
    \mn^{h}_{m} = \binom{m-1}{h-1}    ,\quad 
    \mn^{h}_{n} = \binom{n-1}{h-1}    .
\]
The total number of loops is given by:
\[
    \binom{m-1}{h-1} \binom{n-1}{h-1}    .
\]
Now we can get all sequences by making all $N = m + n$ clockwise
circular permutations on the digits in each loop, but in the process
the $h$ circular permutations displacing an even number of strings
belong to different initial loops.  Then, the total number of sequences 
is:
\begin{equation}\label{e-3}
    \frac{N}{h} \binom{m-1}{h-1} \binom{n-1}{h-1} 
    = \frac{N}{n m} h \binom{m}{h} \binom{n}{h} .
\end{equation}
Finally, taking $2h = \tau$ our formula is:
\begin{equation}\label{e-4}
    \tn^{mn}_{\tau} 
    = \frac{\tau}{2 \mu} \binom{m}{\frac{1}{2} \tau} 
                         \binom{n}{\frac{1}{2} \tau}            ; \quad 
    \tau = 2, 4, 6, \ldots                                      ; \quad 
    \left( \frac{1}{m} + \frac{1}{n} = \frac{1}{\mu} \right)    .
\end{equation}
The even index $\tau$ runs from $2$ to $2 (m \wedge n)$ and
$\tn^{mn}_{\tau}$ is symmetric in the upper indexes $m$, $n$.
\begin{table}[H]\caption{($\tn$-numbers)}\label{t-2}%
\begin{center}%
\setlength{\tabcolsep}{0.25ex}
\begin{tabular}{|r|r|r|r|r|r|r|r|r|r|r|r|r|r|r|r|r|r|r|r|r|r|r|r|r|r|r|r|}\hline
    &        & $N$ &   2 &   3 & \multicolumn{2}{c|}{4} & \multicolumn{2}{c|}{5} & \multicolumn{3}{c|}{6} & \multicolumn{3}{c|}{7} & \multicolumn{4}{c|}{8} & \multicolumn{4}{c|}{9} & \multicolumn{5}{c|}{10} \\ \cline{3-28}
$h$ & $\tau$ & $m$ &\h 1 &\h 1 &\h 2 &\h 1 &\h 2 &\h 1 &\h 3 &\h 2 &\h 1 &\h 3 &\h 2 &\h 1 &\h 4 &\h 3 &\h 2 &\h 1 &\h 4 &\h 3 &\h 2 &\h 1 &\K 5 &\K 4 &\K 3 &\K 2 &\K 1 \\ \hline
  1 &      2 &     &   2 &   3 &   4 &   4 &   5 &   5 &   6 &   6 &   6 &   7 &   7 &   7 &   8 &   8 &   8 &   8 &   9 &   9 &   9 &   9 &  10 &  10 &  10 &  10 &  10 \\ \cline{1-2} \cline{4-28}
  2 &      4 &     &     &     &   2 &     &   5 &     &  12 &   9 &     &  21 &  14 &     &  36 &  32 &  20 &     &  54 &  45 &  27 &     &  80 &  75 &  60 &  35 &     \\ \cline{1-2} \cline{4-28}
  3 &      6 &     &     &     &     &     &     &     &   2 &     &     &   7 &     &     &  24 &  16 &     &     &  54 &  30 &     &     & 120 & 100 &  50 &     &     \\ \cline{1-2} \cline{4-28}
  4 &      8 &     &     &     &     &     &     &     &     &     &     &     &     &     &   2 &     &     &     &   9 &     &     &     &  40 &  25 &     &     &     \\ \cline{1-2} \cline{4-28}
  5 &     10 &     &     &     &     &     &     &     &     &     &     &     &     &     &     &     &     &     &     &     &     &     &   2 &     &     &     &     \\ \hline
\end{tabular}%
\end{center}
\end{table}

From \eqref{e-4} we have the recurrence formula:
\begin{equation}\label{e-5}
    \tn^{mn}_{\tau + 2} 
    = 4 \frac{(m-\frac{1}{2}\tau)(n-\frac{1}{2}\tau)} 
      {\tau (\tau + 2)} \tn^{mn}_{\tau}
\end{equation}
and $\tn^{mn}_{2} =$ number of clockwise permutations of the sequence
\[
    000\stackrel{m}{\cdots}0111\stackrel{n}{\cdots}1 = m + n = N    .
\]
In contrast with further generalizations, \eqref{e-5} shows the 
simplicity of these numbers:
\begin{equation}\label{e-6}
    \begin{split}
        \tn^{mn}_{2} &= N  ; \\
        \tn^{mn}_{4} &= N \cdot \frac{(m-1)(n-1)}{1 \cdot 2} ; \\
        \tn^{mn}_{6} &= N \cdot \frac{(m-1)(n-1)}{1 \cdot 2} %
                          \cdot \frac{(m-2)(n-2)}{2 \cdot 3} ; \\
                     &\;\;\vdots
    \end{split}
\end{equation}
\section{Properties of the Numbers $\tn^{mn}_{\tau}$}\label{s-PNT}

At this point we can relate our problem to the probabilities. Since
$\binom{N}{m}^{-1} \tn^{mn}_{\tau}$ is the probability that a 
sequence have $\tau$ jumps, normalisation gives:
\[
    \sum_{\tau=2}^{2(m \wedge n)} \tn^{mn}_{\tau} = \binom{N}{m}
\]
i.e.:
\begin{equation}\label{e-7}
    \sum_{h=1}^{m \wedge n} h \binom{m}{h} \binom{n}{h} 
    = \frac{m n}{N} \binom{N}{m}
\end{equation}
\begin{equation}\label{e-7p}\tag{\protect\ref{e-7}$'$}
    \sum_{h=1}^{m} h \binom{m}{h}^{2} 
    = \frac{m}{2} \binom{2m}{m}
\end{equation}
which is a well known formula (Knuth \cite[p.\ 59]{Kde1968})%
\footnote{From
    \begin{equation}\label{e-7pp}\tag{\protect\ref{e-7}$''$}
        \sum \binom{r}{h} \binom{s}{m-h} = \binom{r+s}{m}    ,
    \end{equation}
    (Knuth \cite[p.\ 58]{Kde1968}), with $m=r=s$, we have also:
    \begin{equation}\label{e-7ppp}\tag{\protect\ref{e-7}$'''$}
        \sum_{1}^{m} \binom{m}{h}^{2} = \binom{2m}{m}    .
    \end{equation}
    A method which we use later, consists of relating these formulas:
    \[
                \sum_{1}^{m} h \binom{m}{h}^{2} 
              = \sum_{0}^{m} (m-h) \binom{m}{m-h}^{2} 
              = \sum_{0}^{m} (m-h) \binom{m}{h}^{2} 
              = m \sum_{0}^{m} \binom{m}{h}^{2} 
              - \sum_{1}^{m} h \binom{m}{h}^{2}    ,
    \]
    which gives \eqref{e-7ppp} from \eqref{e-7p}.}

For the ordinary moments of the above distribution we have:
\begin{equation}\label{e-8}
    \begin{split}
        \overline{\tau^{r}} &
        = \binom{N}{m}^{-1} \sum \tau^{r} \tn^{mn}_{\tau}    ; \\
        \overline{\tau^{r-1}} &
        = \binom{N}{m}^{-1} \sum (2h)^{r-1} \frac{N}{mn} h 
          \binom{m}{h} \binom{n}{h} 
        = \binom{N}{m}^{-1} \frac{N}{mn} 2^{r-1} 
          \sum^{m \wedge n} h^{r} \binom{m}{h} \binom{n}{h}    .
    \end{split}
\end{equation}
Obviously the probability of jumping is $p = \frac{2mn}{N(N-1)}$, 
(probabilities of finding $(00)$ and $(11)$ are 
$\frac{m(m-1)}{N(N-1)}$ and $\frac{n(n-1)}{N(N-1)}$ respectively); 
then, according to binomial distribution, the probability of jumping 
exactly $\tau$ times is:
\begin{equation}\label{e-9}
    P_{\tau} = 2 \binom{N}{\tau} p^{\tau} (1-p)^{N-\tau} \qquad 
               (\approx \binom{N}{m}^{-1} \tn^{mn}_{\tau})
\end{equation}
(the factor 2 is due to the fact that jumps came by pairs), which shows how 
the theory of binary sequences would be complicated, even for physical 
problems, if one follows statistical formulations, because the exact formulas 
are much more simple than the approximate ones. As is shown in the following 
table, statistical approximation is good enough anyway:
\begin{table}[H]\caption{}\label{t-3}%
\begin{center}%
\newcolumntype{.}{D{.}{.}{2}}
\begin{tabular}{r . . . . . .}\hline
\tstrut $\tau$                         & \K 0 & \K 2  & \K 4  & \K 6   & \K 8  & \h 10 \\ \hline
\tstrut $\binom{10}{5}\tn^{55}_{\tau}$ & 0    & 10    & 80    & 120    & 40    &  2    \\ \hline
\tstrut Binomial                       & 0.15 & 10.66 & 77.71 & 121.42 & 40.65 &  1.46 \\ \hline
\end{tabular}%
\end{center}
\end{table}
\noindent Since both distributions are so close (at least in this example), 
we can expect a very good approximations identifying its moments. Call
$\overline{k^{r}}$ and $\overline{(k)}_{r}$ the ordinary and factorial
moments respectively. We have for a distribution 
$p_{0}, p_{1},\ldots, p_{2m \wedge 2n}$:
\begin{equation}\label{e-10}
    \overline{k^{r}} 
    = \sum_{\ell=0}^{2m \wedge 2n \wedge r} \sn^{\ell}_{r}\, 
                                            \overline{(k)}_{\ell} 
    = \sum_{\ell=0}^{2m \wedge 2n \wedge r} \sn^{\ell}_{r}\, 
                                            G^{(\ell)}(1) 
\end{equation}
the $\sn^{\ell}_{r}$ being the Stirling numbers of second kind and 
$G^{(\ell)}(1)$ the $\ell$ derivative of the generating function $G(t)$ at 
$t=1$; (here we will be interested only in the case 
$2m \wedge 2n \wedge r = r$). The generating function for Bernouilli's 
distribution is $G(t) = (q + pt)^{N}$, so we have for \eqref{e-10}:
\begin{equation}\label{e-11}
    \overline{k^{r}} 
    = \sum_{\ell=0}^{r} \sn^{\ell}_{r} \, (N)_{\ell} \, p^{\ell}
\end{equation}
then, for $r > 1$, (for $r = 1$ we have the formula \eqref{e-7}):
\begin{equation}\label{e-12}
    \sum_{h=1}^{m \wedge n} h^{r} \binom{m}{h} \binom{n}{h} 
    \approx \binom{N}{m \wedge n}
            \frac{ m^{2} n^{2} }{ 2^{r-2} N (N-1) } 
            \sum_{\ell=0}^{m \wedge n} \sn^{\ell}_{r-1} A^{\ell}_{mn}
\end{equation}
with
\[
    A^{1}_{mn} = 1    , \quad
    A^{\ell}_{mn} = \left( \frac{2mn}{N} \right)^{\ell-1}
                    \frac{(N-2)_{\ell-2}}{(N-1)^{\ell-2}}    ; \quad
    (N)_{0}=1    .
\]
As we will see in the sequel, approximation \eqref{e-12} is
extraordinary good for $m=n$, which gives:
\begin{equation}\label{e-13}
    \begin{split}
        \sum_{h=0}^{m} h^{r} \binom{m}{h}^{2} %
            & \approx \binom{2m}{m}
                      \frac{ m^{3} }{ 2^{r-1}(2m-1) } %
                      \left( \sn^{1}_{r-1} + m \sn^{2}_{r-1} %
                             + m^{2} \sn^{3}_{r-1} \frac{2m-2}{2m-1} \right. \\
            & \quad   \left. + \cdots + m^{r-2} \sn^{r-1}_{r-1} %
                             \frac{2m-2}{2m-1} \cdot \: \cdots \: \cdot %
                             \frac{2m-r+2}{2m-1} \right). %
    \end{split}
\end{equation}
For discussion, we write explicitely \eqref{e-13} from $r=0$ to $r=5$:
\begin{align}
    \binom{2m}{m}^{-1} \sum_{h=0}^{m} h^{0} \binom{m}{h}^{2} 
        &\approx 1 
                            \label{e-13,0}\tag{\protect\ref{e-13}$,0$}    \\
    \binom{2m}{m}^{-1} \sum_{h=0}^{m} h^{1} \binom{m}{h}^{2} 
        &\approx \frac{m}{2} 
                            \label{e-13,1}\tag{\protect\ref{e-13}$,1$}    \\
    \binom{2m}{m}^{-1} \sum_{h=0}^{m} h^{2} \binom{m}{h}^{2} 
        &\approx \frac{m^{3}}{2(2m-1)} 
                            \label{e-13,2}\tag{\protect\ref{e-13}$,2$}    \\
    \binom{2m}{m}^{-1} \sum_{h=0}^{m} h^{3} \binom{m}{h}^{2} 
        &\approx \frac{m^{3}}{4(2m-1)}(1+m) 
                            \label{e-13,3}\tag{\protect\ref{e-13}$,3$}    \\
    \binom{2m}{m}^{-1} \sum_{h=0}^{m} h^{4} \binom{m}{h}^{2} 
        &\approx \frac{m^{3}}{8(2m-1)} \left( 1 + 3m + m^{2} 
            \frac{2m-2}{2m-1} \right) 
                            \label{e-13,4}\tag{\protect\ref{e-13}$,4$}    \\
    \binom{2m}{m}^{-1} \sum_{h=0}^{m} h^{5} \binom{m}{h}^{2} 
        &\approx \frac{m^{3}}{16(2m-1)} \left( 1 + 7m 
            + 6m^{2} \frac{2m-2}{2m-1} 
            + m^{3} \frac{(2m-2)(2m-3)}{(2m-1)^{2}} \right) 
                            \label{e-13,5}\tag{\protect\ref{e-13}$,5$}
\end{align}
We know that formulas \eqref{e-13,0} and \eqref{e-13,1} are exact 
(\eqref{e-7ppp}, \eqref{e-7p}) and we will show that \eqref{e-13,2} and 
\eqref{e-13,3} are exact too and also why identity fails after $r=3$.

\begin{proof}
\[
    \begin{split}
            \sum_{1}^{m} h^{2} \binom{m}{h}^{2} 
        &=  \sum_{1}^{m} \left( h \binom{m}{h} \right)^{2} 
            = m^{2} \sum_{1}^{m} \binom{m-1}{h-1}^{2} 
            = m^{2} \binom{2m-2}{m-1} \\
        &=  m^{2} \frac{m(2m-m)}{2m(2m-1)} \binom{2m}{m} 
            = \frac{m^{3}}{2(2m-1)} \binom{2m}{m}    .
    \end{split}
\]
In the general case:
\[
    \sum_{1}^{m} h^{r} \binom{m}{h}^{2} 
    = \sum_{0}^{m} (m-h)^{r} \binom{m}{h}^{2} 
    = \sum_{k=0}^{r} (-)^{k} \binom{r}{k} m^{r-k} 
      \left( \sum_{h=0}^{m} h^{k} \binom{m}{h}^{2} \right)
\]
then, we have $\sum h^{r} \binom{m}{h}^{2}$ by recurrence 
\emph{only for odd values of $r$}, because when $k=r=2$ the last term cancels 
with the left hand side, giving only an identity between moments of lower 
orders. Only the case $r=2$ has escaped from that inconvenient (due to the 
square in $\binom{m}{h}$) and provides for \eqref{e-13,3}. Thus, for $r=3$, 
we have:
\begin{multline*}
    \binom{2m}{m}^{-1} \sum_{0}^{m} h^{3} \binom{m}{h}^{2} \\
    = \binom{2m}{m}^{-1} %
      \left( m^{3} \sum_{0}^{m} \binom{m}{h}^{2} 
           - 3 m^{2} \sum_{1}^{m} h \binom{m}{h}^{2} 
           + 3 m \sum_{1}^{m} h^{2} \binom{m}{h}^{2} 
           - \sum_{0}^{m} h^{3} \binom{m}{h}^{2} \right)
\end{multline*}
i.e.:
\[
    \binom{2m}{m}^{-1} \sum_{0}^{m} h^{3} \binom{m}{h}^{2} 
    = \frac{1}{2} 
      \left( m^{3} - 3 m^{2} \frac{m}{2} + 3 m \frac{m^{3}}{2(2m-1)} 
      \right) 
    = \frac{m^{3} (m+1)}{4 (2m-1)} .
\]
Recurrence fails for $\sum h^{4} \binom{m}{h}^{2}$ and 
so for higher orders.
\end{proof}

Formula \eqref{e-13,4} gives the following approximations. 
For $m=2$: $30/9 \approx 29/9$; for $r=3$: $11.70 \approx 11.61$; 
for $m=10$: $827.40 \approx 827.22$; and 
for $m=15$; $3829.74 \approx 3829.48$. 
Formula \eqref{e-13,5} gives, for $m=10$: $4895.51 \approx 4891.65$.  
Then, the distributions $\binom{2m}{n}^{-1}\tn^{mn}_{\tau}$ and
$2 \binom{2m}{\tau} p^{\tau} (1-p)^{N-\tau}$ have the same variance; higher 
order moments are almost equal in both distributions.

Let us look for these relations in case $m \neq n$: Supose $m < n$ 
(the case $m > n$ is the same); according to \eqref{e-7pp} we have:
\[
      \sum_{h=0}^{m} \binom{m}{h} \binom{n}{h} 
    = \sum_{0}^{m} \binom{m}{h} \binom{n}{n-h} 
    = \binom{m+n}{n} ,
\]
then:
\[
    \begin{split}
        \sum_{h=1}^{m \wedge n} h^{2} \binom{m}{h} \binom{n}{h} 
        &= \sum_{1}^{m} \left( h \binom{m}{h} \right) 
           \left( h \binom{n}{h} \right) 
         = m n \sum_{1}^{m-1} \binom{m-1}{h-1} \binom{n-1}{h-1} \\
        &= m n \sum_{0}^{m-1} \binom{m-1}{h} \binom{n-1}{h} 
         = m n \binom{m+n-2}{m-1} 
         = \frac{m^{2} n^{2}}{N(N-1)} \binom{N}{m}
    \end{split}
\]
which is the formula \eqref{e-12} for $r=2$. For $r=3$, formula \eqref{e-12} 
gives:
\[
    \sum_{h=1}^{m \wedge n} h^{3} \binom{m}{h} \binom{n}{h} 
    \approx \binom{N}{m \wedge n} 
            \left( \frac{1}{2} \frac{m^{2}n^{2}}{N(N-1)} 
                 + \frac{m^{3}n^{3}}{N^{2}(N-1)} \right)
\]
which is only approximate; for $m=4$, $n=2$, ($N=6$) gives 
$56 \approx 58.67$. Of course the approximation will be much better for 
higher $N$. Then, for $m \neq n$, only the mean values are equal in the
above distributions, though higher order moments are very close in
both distributions.

Let us consider the set of all $2^{N}$ sequences of length $N$. One
may expect here complete agreement with probabilities, i.e.\ jumps must
follow the binomial distribution with $p = \frac{1}{2}$. The next
proof shows that this is the case:
\begin{proof}%
From
\[
    \sum_{0 \leq c \leq a} \binom{a-c}{j} \binom{b+c}{k}
    = \binom{a+b+1}{j+k+1}
\]
(Knuth \cite[p.\ 58]{Kde1968}), we have ($b=0$, $j=k$):
\[
    \sum_{0 \leq c \leq a} \binom{a-c}{k} \binom{c}{k} 
    = \binom{a+1}{2k+1}
\]
Then the exact number of sequences with $\tau$ jumps ((01) or (10)) in all 
$2^{N}$ sequences is, according formula \eqref{e-3}:
\begin{equation}\label{e-14}
    \begin{split}
        \sum_{n=0}^{N} \tn^{N-n,n}_{\tau} 
        &= \frac{N}{h} \sum_{n} \binom{(N-2)-(n-1)}{h-1} 
           \binom{n-1}{h-1} \\
        &= \frac{N}{h} \binom{N-1}{2h-1} 
         = 2 \binom{N}{2h} 
         = 2 \binom{N}{\tau} .
    \end{split}
\end{equation}
So, the probability of jumping $\tau$ times: 
$\sum_{n} \tn^{mn}_{\tau} / 2^{N} = \binom{N}{\tau} / 2^{N-1}$ equals that 
given by Bernouilli's distribution \eqref{e-9} with $p = \frac{1}{2}$: 
$2\binom{N}{\tau} \left(\frac{1}{2}\right)^{\tau} %
                  \left(\frac{1}{2}\right)^{N-\tau}$, %
($\tau = 0,2,\ldots,N\ (\mbox{or } N-1)$).
\end{proof}
\subsection*{Application to the Ising model of ferromagnetism}\label{ss-AIMFM}

Let us consider a one-dimensional ferromagnet in absence of external field and 
put $\nu = J/kT$ ($J=$ interaccion energy of antiparallel spins ($-J$ for 
parallel); $k$ is the Boltzmann constant and $T$ the absolute temperature). The 
energy of a given distribution of $N$ spins is 
$J(-(N-\tau)+\tau) = J(-N+2\tau)$, because there are $\tau$ couplings between 
antiparallel spins. The partition function for a fixed number of $N-n$ and $n$ 
spins up and spins down respectively, is:
\[
    Z = \sum_{\tau} \tn^{N-n,n}_{\tau} e^{-(-N+2\tau)\nu}
\]
(this is also the partition function corresponding to binary alloys (ions A 
and B instead of spins up and down)). Fortunately (for the Physicists) one 
must take the sum of all configurations which gives:
\[
    \begin{split}
        Z &= \sum_{\tau} \sum_{n} \tn^{N-n,n}_{\tau} e^{-(-N+2\tau)\nu}
           = 2 e^{N\nu} \sum_{\tau} \binom{N}{\tau} (e^{-2\nu})^{\tau} \\
          &= e^{N\nu} ( (1+e^{-2\nu})^{N} + (1-e^{-2\nu})^{N} )
           = (2 \cosh \nu)^{N} + (2 \sinh \nu)^{N} .
    \end{split}
\]
One wonders why some authors give the solution which contains only the 
hyperbolic cosinus (see, for example, \cite[ch.\ 5]{Tcj1972}; $\sinh$ appears 
in more elaborate discussion \cite[ch.\ 4, p.\ 121]{Mew1964}) though both 
solutions are equal for $N \rightarrow \infty$. We think their countings are 
not accurate enough, for they take the following recurrence: 
$Z_{N} = (2 \cosh \nu) Z_{N-1}$ on grounds that a new spin provides an ammount 
of $\pm J$ in the energy, i.e.\ forgetting a border effect. In fact, consider 
three spins: A up, B down and C up and add a new spin D; the old energy is 
$2J-J=J$, and the new energy is $0$ for D up or $4J$ for D down and not 
$J \pm 2J$ given by the above rule.
\section{Asymptotic Expansions}\label{s-AE}

After $N = 8$ we can use the asymptotic approximation for
$\tn^{mn}_{h}$, with the help of Stirling's formula which gives:
\begin{equation}\label{e-15}
    \binom{m}{h} \approx 2^{m+1} 
    \frac{ e^{ -\frac{ (2h-m)^{2} }{2m} } }{\sqrt{2\pi m}} 
\end{equation}
\eqref{e-15} is surprisely accurate even for low values of $h$ ($< m$);
$m=10$, $h=3$ gives: $120 = \binom{10}{3} \approx 116.1$.  
(For $h= \frac{1}{2} m$ we have 
$\sqrt{\pi_{N}} = \frac{2^{N+1}}{\sqrt{2N}} \binom{N}{\frac{1}{2} N}^{-1}$ 
which is, of course, the Wallis' formula because 
$\binom{N+2}{\frac{1}{2} N + 1} = 4 \binom{N}{\frac{1}{2} N} \frac{N+1}{N+2}$ 
gives $\pi_{N+2} = \pi_{N} \frac{N(N+2)}{(N+1)^{2}}$ and since 
$\pi_{2} = \frac{(2^{3})^{2}}{2 \cdot 2} \binom{2}{1}^{-2} = 4$, we have:
\[
    \pi_{4} = 2 \cdot \frac{2 \cdot 2}{1 \cdot 3} \cdot \frac{4}{3} ; 
    \pi_{6} = 2 \cdot \frac{2 \cdot 2}{1 \cdot 3} 
                \cdot \frac{4 \cdot 4}{3 \cdot 5} \cdot \frac{6}{5} ;
    \ldots ; 
    \pi_{\infty} = \pi = 
        2 \cdot \frac{2 \cdot 2}{1 \cdot 3} 
          \cdot \frac{4 \cdot 4}{3 \cdot 5} 
          \cdot \frac{6 \cdot 6}{5 \cdot 7} \cdot \: \cdots ).
\]
     
We can expect a good approximation for $\tn^{mn}_{h}$ when $N \geq 10$; 
in fact the following formula gives, for $N = 10$ a result very close
to the Table \ref{t-2}:
\[
    \tn^{mn}_{h} 
    \approx h \left( \frac{1}{m} + \frac{1}{n} \right) 
    \frac{2^{m+n+2}}{2 \pi \sqrt{mn}}
    e^{- 2 h^{2} \left( \frac{1}{m} + \frac{1}{n} \right) 
       + 4 h - \frac{1}{2} (m+n)}
\]
or:
\begin{equation}\label{e-16}
\tn^{mn}_{h} \approx \frac{2 N e^{aN}}{\pi (mn)^{3/2}} h 
    e^{-4 \left( \frac{h^{2}}{2 \mu} \right) - h} ; \quad
(a = - \frac{1}{2} + \log 2 \approx 0.193147).
\end{equation}
In terms of $\tau$'s
\begin{equation}\label{e-17}
    \tn^{mn}_{\tau} 
    \approx \tau 
    \frac{ e^{- \frac{\tau^{2}}{2\mu} + 2\tau + aN} }
         { \pi \mu^{3/2} \sqrt{N} } .
\end{equation}
For $N = 10$, $m = n = 5$: $\tn^{55}_{h} 
= 0.3514 \, h \exp \big( - 4 \big( \frac{h^{2}}{5} - h \big) \big)$ 
and we have the following results, as is shown in Figure \ref{f-1}.
\begin{figure}[H]\caption{Asymptotic expansion for $\tn^{55}_{h}$}\label{f-1}
\begin{center}
\setlength{\unitlength}{3ex}
\begin{picture}(18,20)(-2,-2)
%
\put(0,0){\vector(1,0){14}}                   
    \multiput(2,0)(2,0){6}{\line(0,1){0.2}}   
\put(0,0){\vector(0,1){16}}                   
    \multiput(0,1)(0,1){15}{\line(1,0){0.2}}  
    \multiput(0,5)(0,5){3}{\line(1,0){0.4}}   
\put(-1, 0){\makebox(0,0){$\phantom{00}0$}}   
\put(-1, 5){\makebox(0,0){$\phantom{0}50$}}   
\put(-1,10){\makebox(0,0){$100$}}             %
\put(-1,15){\makebox(0,0){$150$}}             %
\put(0,17){\makebox(7.75,0){$\tn^{55}_{x} = 0.3514 \, x %
           \exp \big( - 4 \big( \frac{x^{2}}{5} - x \big) \big) $}}
\put( 2,-0,75){\makebox(0,0){$1$}}   
\put( 4,-0,75){\makebox(0,0){$2$}}   %
\put( 6,-0,75){\makebox(0,0){$3$}}   %
\put( 8,-0,75){\makebox(0,0){$4$}}   %
\put(10,-0,75){\makebox(0,0){$5$}}   %
\put(12,-0,75){\makebox(0,0){$6$}}   %
\put(14,-0,75){\makebox(0,0){$h$}}   
\put( 2,-1.5){\makebox(0,0){$2$}}    %
\put( 4,-1.5){\makebox(0,0){$4$}}    %
\put( 6,-1.5){\makebox(0,0){$6$}}    %
\put( 8,-1.5){\makebox(0,0){$8$}}    %
\put(10,-1.5){\makebox(0,0){$10$}}   %
\put(12,-1.5){\makebox(0,0){$12$}}   %
\put(14,-1.5){\makebox(0,0){$\tau$}} 
\put( 2, 1){\circle{0.4}} \put(3, 1){\makebox(0,0){$10$}}
\put( 4, 8){\circle{0.4}} \put(5, 8){\makebox(0,0){$80$}}
\put( 6,12){\circle{0.4}} \put(7.2,12){\makebox(0,0){$120$}}
\put( 8, 4){\circle{0.4}} \put(9, 4){\makebox(0,0){$40$}}
\put(10,0.2){\circle{0.4}}\put(10.8,0.5){\makebox(0,0){$2$}}
\put(12,0.0){\circle{0.4}}\put(12.8,0.5){\makebox(0,0){$0$}}
\path%
(0,0)%
(0.1,0.002141717)%
(0.2,0.005200501)%
(0.3,0.009433056)%
(0.4,0.015149)%
(0.5,0.022715)%
(0.6,0.032569)%
(0.7,0.045219)%
(0.8,0.061255)%
(0.9,0.081356)%
(1,0.106292)%
(1.1,0.136935)%
(1.2,0.174254)%
(1.3,0.219326)%
(1.4,0.273326)%
(1.5,0.337531)%
(1.6,0.413309)%
(1.7,0.502111)%
(1.8,0.605455)%
(1.9,0.724909)%
(2,0.862073)%
(2.1,1.018545)%
(2.2,1.195895)%
(2.3,1.395631)%
(2.4,1.619158)%
(2.5,1.867739)%
(2.6,2.142448)%
(2.7,2.44413)%
(2.8,2.773353)%
(2.9,3.130361)%
(3,3.515037)%
(3.1,3.926858)%
(3.2,4.36487)%
(3.3,4.827651)%
(3.4,5.3133)%
(3.5,5.81942)%
(3.6,6.343124)%
(3.7,6.881044)%
(3.8,7.429352)%
(3.9,7.983798)%
(4,8.539752)%
(4.1,9.09227)%
(4.2,9.636155)%
(4.3,10.166039)%
(4.4,10.676469)%
(4.5,11.161999)%
(4.6,11.617283)%
(4.7,12.037179)%
(4.8,12.416839)%
(4.9,12.751805)%
(5,13.038096)%
(5.1,13.272287)%
(5.2,13.451576)%
(5.3,13.57384)%
(5.4,13.63768)%
(5.5,13.642443)%
(5.6,13.588233)%
(5.7,13.475912)%
(5.8,13.307071)%
(5.9,13.083998)%
(6,12.809628)%
(6.1,12.487478)%
(6.2,12.121575)%
(6.3,11.716373)%
(6.4,11.276666)%
(6.5,10.807494)%
(6.6,10.314052)%
(6.7,9.801595)%
(6.8,9.275349)%
(6.9,8.740427)%
(7,8.201752)%
(7.1,7.663987)%
(7.2,7.131479)%
(7.3,6.608204)%
(7.4,6.097736)%
(7.5,5.603216)%
(7.6,5.127334)%
(7.7,4.67233)%
(7.8,4.239992)%
(7.9,3.831669)%
(8,3.448292)%
(8.1,3.090404)%
(8.2,2.758182)%
(8.3,2.451483)%
(8.4,2.169874)%
(8.5,1.912677)%
(8.6,1.679002)%
(8.7,1.467794)%
(8.8,1.277863)%
(8.9,1.107925)%
(9,0.956631)%
(9.1,0.822598)%
(9.2,0.704435)%
(9.3,0.600767)%
(9.4,0.510251)%
(9.5,0.431594)%
(9.6,0.363564)%
(9.7,0.305002)%
(9.8,0.254825)%
(9.9,0.21203)%
(10,0.1757)%
(10.1,0.144999)%
(10.2,0.119174)%
(10.3,0.097547)%
(10.4,0.079519)%
(10.5,0.064558)%
(10.6,0.052198)%
(10.7,0.042032)%
(10.8,0.033708)%
(10.9,0.026922)%
(11,0.021415)%
(11.1,0.016965)%
(11.2,0.013385)%
(11.3,0.010517)%
(11.4,0.008230167)%
(11.5,0.006414365)%
(11.6,0.004978845)%
(11.7,0.003848878)%
(11.8,0.002963267)%
(11.9,0.002272161)%
(12,0.001735160)
\end{picture}
\end{center}
\end{figure}
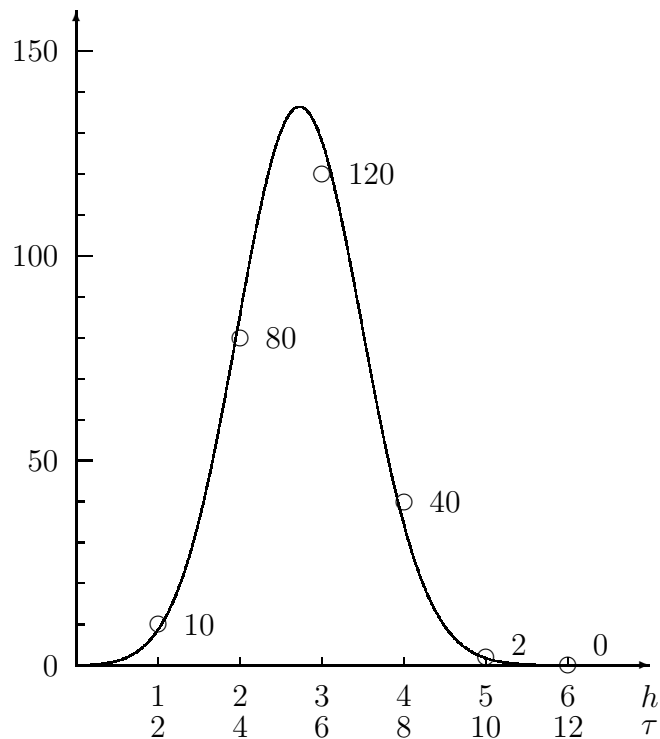
\begin{table}[H]\caption{}\label{t-4}%
\begin{center}%
\newcolumntype{.}{D{.}{.}{2}}%
\begin{tabular}{r . . . . . .}\hline
\tstrut $\tau$            & \K 2 & \K  4 & \K   6 & \K  8 & \h 10 & \h 12 \\ \hline
\tstrut $\tn^{55}_{\tau}$ &   10 &    80 &    120 &    40 &     2 &     0 \\ \hline
\tstrut Asymptotic        & 8.62 & 80.40 & 127.95 & 34.44 &  1.76 &  0.02 \\ \hline
\end{tabular}%
\end{center}%
\end{table}%

At the origin, the derivative of \eqref{e-16} is not zero but 
$2 (\pi \mu \sqrt{mn})^{-1} e^{aN}$ and the curve takes, though very small, 
negative values for $x < 0$.

On considering the set of all $2^{N}$ sequences, the distribution of
jumps, for a large $N$ is, with \eqref{e-14} and \eqref{e-15}:
\begin{equation}\label{e-18}
    \sum_{n} \tn^{N-n,n}_{\tau} 
    = 2 \binom{N}{\tau} 
    \approx \frac{2^{N+2}}{\sqrt{2 \pi N}} 
            e^{ - \frac{(2\tau - N)^{2}}{2N} } .
\end{equation}
Here the gaussian has nothing to do with statistical distributions,
because in this context \eqref{e-18} gives fluctuations, respect to 
the exact number of jumps, due only to the fact that $N$ is not large 
enough. Of course \eqref{e-18} becomes a normal probability 
distribution with standard error $\pm\sqrt{N}$ when the sequences of 
$N$ digits are taken from a table of binary random digits.
\section{The Numbers $\tn^{mn}_{h}(00)$ and $\tn^{mn}_{h}(11)$.\\
         Generalized Fibonacci Numbers}\label{s-TNTTGFN}

Let $U$ be an ordered sequence of $r$ zeros and $s$ ones, 
$|U| = r + s < N$; being $\tau$ the amount of $(01)$'s or $(10)$'s in 
the sequence $S_{i}$, our general problem will consist in finding the 
number of times, $h = h(U)$, a fixed string $U$ is contained in $S_{i}$ 
in the following sense: If 
$S_{i} = \alpha_{1}\; \alpha_{2}\; \alpha_{3}\; \cdots\; \alpha_{N}$, 
($\alpha_{j} = 0, 1$), take the $N$ strings:
\begin{multline*}
    (\alpha_{1}\; \alpha_{2}\; \alpha_{3}\; \cdots\; \alpha_{r+s} ) , 
    (\alpha_{2}\; \alpha_{3}\; \alpha_{4}\; \cdots\; \alpha_{r+s+1} ) , 
    \ldots , \\
    (\alpha_{N-1}\; \alpha_{N}\; \alpha_{1}\; \cdots\; \alpha_{r+s-2} ) , 
    (\alpha_{N}\; \alpha_{1}\; \alpha_{2}\; \cdots\; \alpha_{r+s-1} ) , 
\end{multline*}
then $h$ is the number of these strings having the form $U$. For instance, 
take in \eqref{e-2} $U = (1100)$, then $h = 3$. We call $\tn^{mn}_{h}U$ the 
number of sequences which have $h$ times the string $U$.

Now, since after any $(01)$ in a sequence $S$ sooner or later another 
$0$ appears (Ising model!), then $h(10) = h(01) = \frac{1}{2} \tau$, i.e., 
\begin{equation}\label{e-19}
    \tn^{mn}_{h}(01) = \tn^{mn}_{h}(10) = \tn^{mn}_{\tau} .
\end{equation}

The numbers $\tn^{mn}_{h}(00)$ and $\tn^{mn}_{h}(11)$ are found by a 
simple consideration: From all $N$ pairs succesively taken in $S$, 
$m =$ the number of zero digits $=$ the number $h + \frac{1}{2} \tau$, 
of them having the form $(00)$ or $(01)$, so the $\tau$ corresponding 
to $S$ is $2(m-h)$, and we have:
\begin{equation}\label{e-20}
    \tn^{mn}_{h}(00) = \tn^{mn}_{2(m-h)} ,\qquad\text{similarly:}\qquad
    \tn^{mn}_{h}(11) = \tn^{mn}_{2(n-h)} 
\end{equation}
then:
\begin{equation}\label{e-21}
    \tn^{mn}_{h}(00) = N \frac{m-h}{mn} \binom{m}{h} \binom{n}{m-h} ,\qquad 
    \tn^{mn}_{h}(11) = N \frac{n-h}{mn} \binom{m}{m-h} \binom{n}{h}
\end{equation}
which gives the obvious identity: 
$\tn^{mn}_{h}(00) = \tn^{mn}_{h}(11)$.

\begin{normalization}
Since the only values of $h$ for the 
$\tn^{mn}_{h}(00)$ are between $m$ and $m - (m \wedge n)$, we have, 
(the same for $\tn^{mn}_{h}(11)$):
\[
    \sum_{m-h=0}^{m \wedge n} \tn^{mn}_{h}(00) 
    = \sum_{m-h=0}^{m \wedge n} \tn^{mn}_{2(m-h)} 
    = \sum_{j=0}^{m \wedge n} \tn^{mn}_{2j}
    = \sum_{\frac{1}{2}\tau=0}^{m \wedge n} \tn^{mn}_{\tau}
    = \binom{N}{m} .
\]
\end{normalization}

Let $\gfn(N;r,h)$ be the number of subsets of $X = \{ 1,2,3,\ldots,N \}$ 
that do contain exactly $h$ times $r$ consecutive integers modulo $N$; as in 
Berge \cite[p.~31]{Bc1971} let us associate with each subset $Y \subseteq X$ 
a ``word'' $S = \alpha_{1} , \alpha_{2} , \alpha_{3} , \ldots , \alpha_{N}$, 
where $\alpha_{i} = 0$ if $i \not\in Y$ and $\alpha_{i} = 1$ if 
$i \in Y$. The mapping between subsets $Y$ and words $S$ being 
bijective, instead of counting subsets, we can count the number of the 
sequences $S$ in which the string $(11\stackrel{r}{\cdots}1)$ is 
repeated exactly $h$ times. We call \emph{Generalized corrected 
Fibonacci numbers} the number of these subsets, i.e.:
\begin{equation}\label{e-22}
    \gfn(N;r,h) 
    = \sum_{n} \tn^{N-n,n}_{h}(11\stackrel{r}{\cdots}1)
\end{equation}
then:
\begin{equation}\label{e-23}
    \gfn(N;2,h) 
    = \sum_{n} \tn^{N-n,n}_{h}(11) 
    = \sum_{n=h+1}^{\left[ \frac{N+h}{2} \right]} 
      \frac{(n-h)N}{(N-n)n} \binom{N-n}{N-h} \binom{n}{h} .
\end{equation}
\begin{example}
$X = \{ 1,2,3,4 \}$, $r = 2$, $h = 2$. In this case, in \eqref{e-23} the only 
index is $n = 3$, then $\gfn(4;2,2) = 4$. The subsets are: $\{1,2,3\}$, 
$\{2,3,4\}$, $\{3,4,1\}$, $\{2,4,1\}$. For $h = 0$ one obtains the corrected 
Fibonacci numbers: (see Berge, \cite[p.~32]{Bc1971})
\[
    \gfn(N;2,0) 
    = \gfn_{N} 
    = \sum_{n=1}^{\left[ \frac{1}{2} N \right]} 
      \frac{N}{N-n} \binom{N-n}{n}.
\]
\end{example}

Next, to each of $\binom{N}{m} = \binom{N}{n}$ sequences we will 
associate a \emph{type}, which is the pair of partitions
\[
    ( \beta^{\alpha_{1}}_{1}, \beta^{\alpha_{2}}_{2}, \ldots, 
      \beta^{\alpha_{s}}_{s} ) , \qquad
    ( \beta'_{1}{}^{\alpha'_{1}}, \beta'_{2}{}^{\alpha'_{2}}, \ldots, 
      \beta_{s}'{}^{\alpha'_{s}} ) , 
\]
with height $h$, i.e.: 
$\sum_{1}^{s} \alpha_{i} = \sum_{1}^{s'} \alpha'_{i} = h$, corresponding to the 
tableaux of $0$'s and $1$'s, respectively, disregarding the order of the 
blocks, i.e.: in descending order of the length of blocks, as is customary. 
Two types are equal iff the have both Young's tableaux equal. For instance, 
for $m=4$, $n=3$, the $35$ sequences of the Table \ref{t-1} have the following 
four types:
\[
    \mathbf{1} = \{ (4),(3) \} , \quad 
    \mathbf{2} = \{ (3,1),(2,1) \} , \quad 
    \mathbf{3} = \{ (2^{2}),(2,1) \} , \quad 
    \mathbf{4} = \{ (2,1^{2}),(1^{3}) \} . 
\]
Thus,
\[
    \type \{ 0011100 \} = \mathbf{1} , \quad 
    \type \{ 0110100 \} = \type \{ 0001011 \} = \mathbf{2} , \quad 
    \text{etc.}
\]
In total, for a fixed $m$, $n$, there are the number
\begin{equation}\label{e-24}
    \sum_{h=1}^{m \wedge n} \pn^{h}_{m}\pn^{h}_{n} 
\end{equation}
of types. The great number of sequences compared to the small number of 
types is due to the $\mn$-partitions; in fact, the number of sequences 
per type is:
\begin{equation}\label{e-25}
    \frac{N\, h!\, (h-1)!}
         { \alpha_{1}!\, \alpha_{2}! \cdots \alpha_{s}!\, 
           \alpha'_{1}!\, \alpha'_{2}! \cdots \alpha'_{s}! }
\end{equation}
because the $\mn$-partitions are found making the sum: 
$h! \sum \frac{1}{\alpha_{1}!\, \alpha_{2}! \cdots \alpha_{s}!}$ of all 
partitions with same $h$. This clarifies the meaning of the formulas:
\begin{equation}\label{e-26}
    \begin{split}
        \tn^{mn}_{h}(01) &= \tn^{mn}_{h}(10) 
                          = \frac{N}{n} \mn^{h}_{m} \mn^{h}_{n} ,     \\
        \tn^{mn}_{h}(00) &= \frac{N}{n} \mn^{m-h}_{m} \mn^{m-h}_{n} , \\
        \tn^{mn}_{h}(11) &= \frac{N}{n} \mn^{n-h}_{m} \mn^{n-h}_{n} .
    \end{split}
\end{equation}
\section{The Numbers $\tn^{mn}_{h}U$ for $|U|=3$}\label{s-TNTU3}

The number of formulas can be reduced by some general considerations; 
for example, if $U'$ is the string obtained by interchanging in $U$ the 
$0$'s and the $1$'s we have: $\tn^{mn}_{h}U' = \tn^{nm}_{h}U$. 
Furthermore, since sequences are unbounded, after any string 
${(011\stackrel{k}{\cdots}1)}$ in the sequence sooner or later another 
zero appears, then 
$\tn^{mn}_{h}{(11\stackrel{k}{\cdots}10)} 
= \tn^{mn}_{h}{(011\stackrel{k}{\cdots}1)}$; also, 
$\tn^{mn}_{h}{(00\stackrel{k}{\cdots}01)} 
= \tn^{mn}_{h}{(10\stackrel{k}{\cdots}00)}$, so
$\tn^{mn}_{h}{(11\stackrel{k}{\cdots}10)} 
= \tn^{nm}_{h}{(10\stackrel{k}{\cdots}00)}$, etc. For $U=3$ we have 
the following identities
\[
    \begin{split}
        \tn^{mn}_{h}(001) &= \tn^{mn}_{h}(100) 
                           = \tn^{nm}_{h}(110) 
                           = \tn^{nm}_{h}(011) ; \\
        \tn^{mn}_{h}(101) &= \tn^{nm}_{h}(010) ; \\
        \tn^{mn}_{h}(111) &= \tn^{nm}_{h}(000) .
    \end{split}
\]
Thus we only need formulas for the three expressions:%
\begin{equation}\label{e-27}
    \tn^{mn}_{h}(001),\quad \tn^{mn}_{h}(101),\quad \tn^{mn}_{h}(000).
\end{equation}

Let $\pn^{h}_{m}$be the number of partitions, disregarding order, of 
weight $m$ and dimension $h$; as is well known, the process of deleting 
the first column on all associated Young's Tableaux gives:
\begin{equation}\label{e-28}
    \pn^{h}_{m} = \sum_{\ell=1}^{s} \pn^{\ell}_{m-h} \qquad (s=h\wedge(m-h)) .
\end{equation}
Let us try a similar expansion for $\mn$-partitions:
\begin{equation}\label{e-29}
    \mn^{h}_{m} = \sum_{\ell=1}^{s} a_{\ell} \mn^{\ell}_{m-h} ,
\end{equation}
formula \eqref{e-7} gives:
\begin{equation}\label{e-30}
    \sum_{\ell=1}^{s} \binom{h}{\ell} \mn^{\ell}_{m-h} 
    = \frac{1}{m-h} \sum_{\ell=1}^{s} \ell \binom{h}{\ell} \binom{m-h}{\ell}
    = \frac{h}{m} \binom{m}{h} 
    = \mn^{h}_{m}  .
\end{equation}
Then the coefficients $a_{\ell}$ in the formal expansion \eqref{e-29} 
are equal to $\binom{h}{\ell}$. We shall prove that this expansion is 
not only formal which proves the following lemma: Define 
$c^{i}_{jk} = \frac{k}{i-j} \binom{j}{k} \binom{i-j}{k}$ for 
$i \neq j$ and $c^{i}_{ik} = i \delta_{0k}$, so that \eqref{e-30} has 
the form: $\sum c^{i}_{jk} = \mn^{j}_{i}$. We have,
\begin{lemmaoneone}
The number of $\mn$-partitions of $m$ in $h$ parts 
is $\mn^{h}_{m}$ but after deleting the first column in all 
Young-De Moivre Tableaux we have $\mn$-partitions of dimensions 
$\ell = 0,1,2,\ldots,h\wedge(m-h)$, which number is $c^{m}_{h\ell}$.
\end{lemmaoneone}
\noindent In order that the counting makes sense we must define 
$c^{m}_{m0} = 1$. Further, it is important to take in consideration the 
fact that deleting the first column is made for counting purposes only. 
Consider, for instance, this process for $m=6$, $h=2$; we have the 
$\binom{5}{1}$ following tableaux:
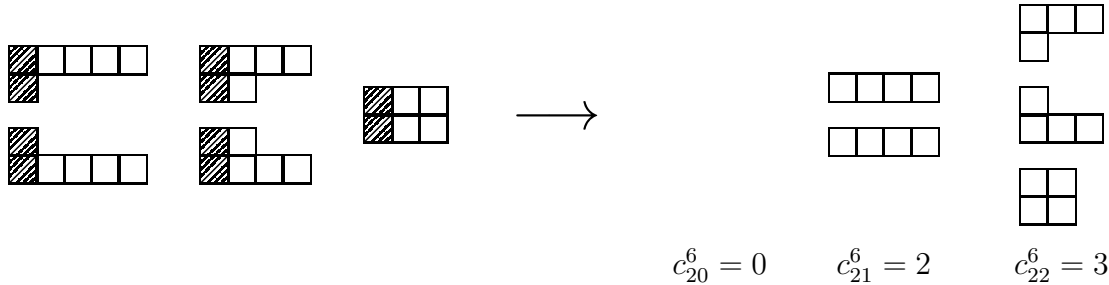
\setlength{\unitlength}{2ex}
\newsavebox{\shbx}
\savebox{\shbx}(1,1){\thinlines%
                     \drawline(0,0.75)(0.25,1)%
                     \drawline(0,0.50)(0.50,1)%
                     \drawline(0,0.25)(0.75,1)%
		     \drawline(0,0)(1,1)%
 		     \drawline(0.25,0)(1,0.75)%
		     \drawline(0.50,0)(1,0.50)%
                     \drawline(0.75,0)(1,0.25)}
\begin{center}
\begin{picture}(40,10)(0,-2)
%
\multiput(0,5.5)(1,0){5}{\framebox(1,1){}}\multiput(0,5.5)(1,0){1}{\usebox{\shbx}}
\multiput(0,4.5)(1,0){1}{\framebox(1,1){}}\multiput(0,4.5)(1,0){1}{\usebox{\shbx}}
\multiput(0,2.5)(1,0){1}{\framebox(1,1){}}\multiput(0,2.5)(1,0){1}{\usebox{\shbx}}
\multiput(0,1.5)(1,0){5}{\framebox(1,1){}}\multiput(0,1.5)(1,0){1}{\usebox{\shbx}}
\multiput(7,5.5)(1,0){4}{\framebox(1,1){}}\multiput(7,5.5)(1,0){1}{\usebox{\shbx}}
\multiput(7,4.5)(1,0){2}{\framebox(1,1){}}\multiput(7,4.5)(1,0){1}{\usebox{\shbx}}
\multiput(7,2.5)(1,0){2}{\framebox(1,1){}}\multiput(7,2.5)(1,0){1}{\usebox{\shbx}}
\multiput(7,1.5)(1,0){4}{\framebox(1,1){}}\multiput(7,1.5)(1,0){1}{\usebox{\shbx}}
\multiput(13,4)(1,0){3}{\framebox(1,1){}}\multiput(13,4)(1,0){1}{\usebox{\shbx}}
\multiput(13,3)(1,0){3}{\framebox(1,1){}}\multiput(13,3)(1,0){1}{\usebox{\shbx}}
\put(20,3.85){\makebox(0,0){\LARGE$\longrightarrow$}}
\put(24,-2){\makebox(4,1){$c^{6}_{20} = 0$}}
\multiput(30,4.5)(1,0){4}{\framebox(1,1){}}
\multiput(30,2.5)(1,0){4}{\framebox(1,1){}}
\put(30,-2){\makebox(4,1){$c^{6}_{21} = 2$}}
\multiput(37,7)(1,0){3}{\framebox(1,1){}}
\multiput(37,6)(1,0){1}{\framebox(1,1){}}
\multiput(37,4)(1,0){1}{\framebox(1,1){}}
\multiput(37,3)(1,0){3}{\framebox(1,1){}}
\multiput(37,1)(1,0){2}{\framebox(1,1){}}
\multiput(37,0)(1,0){2}{\framebox(1,1){}}
\put(37,-2){\makebox(3,1){$c^{6}_{22} = 3$}}
\end{picture}
\end{center}
then $c^{6}_{21} = 2$ means that the tableau%
\setlength{\unitlength}{2ex}
\begin{picture}(4,1)
\multiput(0,0)(1,0){4}{\framebox(1,1){}}
\end{picture}
shall be counted twice because it has different sites with respect to 
the first column.
\begin{proof}[Proof of Lemma I${}^{(1)}$]
Before deleting the first 
column, the number of tableaux which will give tableaux of dimension 
$\ell$ after the cut, is (Figure~\ref{f-2}):%
\begin{figure}[H]\caption{}\label{f-2}
\begin{center}
\setlength{\unitlength}{2.25ex}
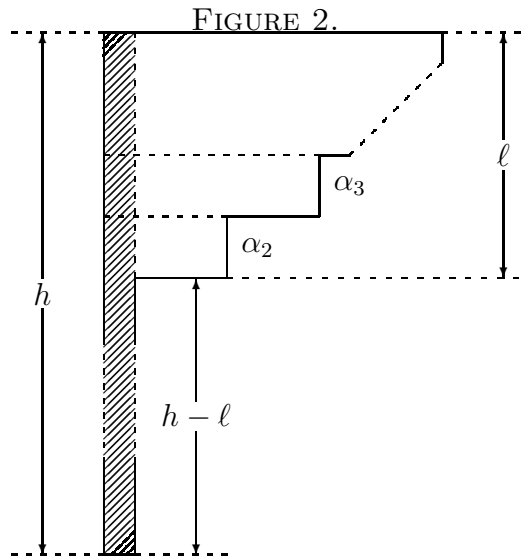
\begin{picture}(17,17)
%
\dashline{0.2}(0,0)(3,0)\drawline(3,0)(4,0)\dashline{0.2}(4,0)(7,0)
\drawline(4,9)(7,9)\dashline{0.2}(7,9)(17,9)
\dashline{0.2}(3,11)(7,11)\drawline(7,11)(10,11)
\dashline{0.2}(3,13)(10,13)\drawline(10,13)(11,13)
\dashline{0.2}(0,17)(3,17)\drawline(3,17)(14,17)\dashline{0.2}(14,17)(17,17)
\drawline(3,0)(3,3)\dashline{0.2}(3,3)(3,7)\drawline(3,7)(3,17)
\drawline(4,0)(4,3)\dashline{0.2}(4,3)(4,7)
\drawline(4,7)(4,9)\dashline{0.2}(4,9)(4,17)
\drawline(7,9)(7,11)
\drawline(10,11)(10,13)
\drawline(14,16)(14,17)
\dashline{0.2}(11,13)(14,16)
\put(1,8.5){\makebox(0,0){$h$}}
\put(1,9){\vector(0,1){8}}\put(1,8){\vector(0,-1){8}}
\put(6,4.5){\makebox(0,0){$h - \ell$}}
\put(6,5){\vector(0,1){4}}\put(6,4){\vector(0,-1){4}}
\put(16,13){\makebox(0,0){$\ell$}}
\put(16,13.5){\vector(0,1){3.5}}\put(16,12.5){\vector(0,-1){3.5}}
\put(8,10){\makebox(0,0){$\alpha_{2}$}}
\put(11,12){\makebox(0,0){$\alpha_{3}$}}
\multiput(3,0)(0,0.25){65}{\drawline(0,0)(1,1)}
\drawline(3.25,0)(4,0.75)    
\drawline(3.50,0)(4,0.50)    
\drawline(3.75,0)(4,0.25)    
\drawline(3,16.25)(3.75,17)  
\drawline(3,16.50)(3.50,17)  
\drawline(3,16.75)(3.25,17)  
\end{picture}%
\end{center}%
\end{figure}%
\[
    \mathcal{N} 
    = \sum \frac{h!}{(h-\ell)!\, \alpha_{2}!\, \alpha_{3}! \cdots} 
    = \frac{h!}{(h-\ell)!} \sum \frac{1}{\alpha_{2}! \, \alpha_{3}! \cdots}
\]
where the sum is extended over all tableaux. After deleting the first 
column, these tableaux will give $\mn^{\ell}_{m-h}$ tableaux of 
dimension $\ell$, which number is:
\[
    \ell! \sum \frac{1}{\alpha_{2}! \, \alpha_{3}! \cdots} 
    \quad (= \mn^{\ell}_{m-h})
\]
then:
\[
    \mathcal{N} 
    = \frac{h!}{(h-\ell)!} \frac{\mn^{\ell}_{m-h}}{\ell!} 
    = \binom{h}{\ell} \mn^{\ell}_{m-h} 
    = c^{m}_{h\ell} .\qedhere
\]
\end{proof}
Due to importance of coefficients $c^{i}_{jk}$ we give the tables (for 
matrices of low $k$) in Appendix A. It is immediate to find a 
recurrence relation:
\begin{equation}\label{e-31}
    c^{i}_{j,k+1} = c^{i}_{jk} \frac{(j-k)(i-j-k)}{k(k+1)} .
\end{equation}

Now it is a simple task to find $\tn^{mn}_{h}(001)$; as in the case of 
$\tn^{mn}_{h}(01)$ we match tableaux of zeros with tableaux of ones of 
same height $h$. Since strings $(001)$ are originated from the rows 
of $0$'s tableaux which remain non empty after deleting the first 
column, here in the first formula \eqref{e-26} instead of $\mn^{h}_{m}$ 
we put $c^{m}_{h\ell}$, which gives a number $\ell$ of $(001)$ strings; 
these $\ell$ strings are a subset of the $h$ strings $(01)$, 
($h \geq \ell$), because both $(001)$ and $(101)$ give $(01)$. Then we 
have an intermediate formula for $(001)$:
\begin{untitled}
    The number of sequences in which the string $(01)$ appears $h$ 
    times and the string $(001)$ appears $\ell$ times is given by:
    \begin{equation}\label{e-32}
        \tn^{mn}_{h,\ell}(01;001) = \frac{N}{h} c^{m}_{h\ell} \mn^{h}_{n} 
    \end{equation}
\end{untitled}
\noindent (Note that: $\sum_{\ell} \tn^{mn}_{h \ell}(01;001) 
= \frac{N}{h} \mn^{h}_{n} \sum_{\ell} c^{m}_{h\ell}$, which is formula 
\eqref{e-26}). Thus $\tn^{mn}_{\ell}(001)$ is found by disregarding 
$h$ in \eqref{e-32}, i.e.:
\begin{subequations}
\begin{equation}\label{e-33a}
    \tn^{mn}_{\ell}(001) = N \sum_{h} \frac{c^{m}_{h\ell}}{h} \mn^{h}_{n} 
                         = \frac{N}{h} \sum_{h} c^{m}_{h\ell} \binom{n}{h} .
\end{equation}

Let us find the exact limits on the sum; obviously we have: 
$\ell = 0,1,2,\ldots,\left[ \frac{m}{2} \right] \wedge n$, and the 
inferior limit corresponds to $h = \ell$. Condition $h \leq m,n$ must 
be restricted because in $c^{i}_{jk}$ the low index $k$ runs from $k=0$ 
to $k=i-j$, then $h \leq m-\ell,n$ and finally:
\begin{equation}\label{e-33b}
    \begin{split}
        \tn^{mn}_{\ell}(001)  
        &= \frac{N}{n} \sum_{h=\ell}^{(m-\ell)\wedge n} 
           c^{m}_{h\ell} \binom{n}{h} \\
        &= \frac{N}{n} \sum_{h=\ell}^{(m-\ell)\wedge n} 
           \frac{\ell}{m-h} \binom{h}{\ell} \binom{m-h}{\ell} 
           \binom{n}{h} .
    \end{split}
\end{equation}
In that follows some binomial coefficients will be taken outside the 
sum by help of the formula: 
$\binom{a}{c} \binom{c}{b} = \binom{a}{b} \binom{a-b}{c-b}$, then 
\eqref{e-33b} is:
\begin{equation}\label{e-33c}
    \tn^{mn}_{\ell}(001)  
    = \frac{N}{n} \binom{n}{\ell} \sum_{h=\ell}^{(m-\ell)\wedge n} 
      \binom{m-h-1}{\ell-1} \binom{n-\ell}{h-\ell}  ,
      \qquad \binom{0}{0} = 1  .
\end{equation}
\end{subequations}    

Let us take an example which can be checked by hand on the $70$ sequences 
of $4$ zeros and $4$ ones. From the tables of Appendix A we have:
\vspace{0.5\baselineskip}
\begin{center}
    $c^{4}_{h\ell} =$ 
\begin{tabular}{|r|r|r|r|r|r|}\hline
        & \multicolumn{5}{|c|}{$\ell$} \\ \hline
    $h$ & \h 0 & \h 1 & \h 2 & \h 3 & \h 4 \\ \hline
    0   &   &   &   &   &   \\ \hline
    1   &   & 1 &   &   &   \\ \hline
    2   &   & 2 & 1 &   &   \\ \hline
    3   &   & 3 &   &   &   \\ \hline
    4   & 4 &   &   &   &   \\ \hline
\end{tabular}
\quad then $\tn^{44}_{h\ell}(01;001) 
            = \dfrac{8}{4} c^{4}_{h\ell}\dbinom{4}{h} =$ 
\begin{tabular}{|r|r|r|r|r|r|}\hline
        & \multicolumn{5}{|c|}{$\ell$} \\ \hline
    $h$ & \h 0 & \h 1 & \h 2 & \h 3 & \h 4 \\ \hline
    0   &   &    &    &   &   \\ \hline
    1   &   &  8 &    &   &   \\ \hline
    2   &   & 24 & 12 &   &   \\ \hline
    3   &   & 24 &    &   &   \\ \hline
    4   & 2 &    &    &   &   \\ \hline
\end{tabular}
\end{center}
\vspace{0.5\baselineskip}
The sum rows and columns gives $\tn^{44}_{\ell}(01)$ and 
$\tn^{44}_{\ell}(001)$ respectively:
\vspace{0.5\baselineskip}
\begin{center}
\begin{tabular}{|r|r|r|r|r|r|r|}\hline
        & \multicolumn{5}{|c|}{$\ell$} & \\ \hline
    $h$ & \h 0 & \h 1 & \h 2 & \h 3 & \h 4 & $\tn^{44}_{\ell}(01)$ \\ \hline
    0   &      &      &      &      &      &  0 \\ \hline
    1   &      &    8 &      &      &      &  8 \\ \hline
    2   &      &   24 &   12 &      &      & 36 \\ \hline
    3   &      &   24 &      &      &      & 24 \\ \hline
    4   &    2 &      &      &      &      &  2 \\ \hline
    $\tn^{44}_{\ell}(001)$ & 2 & 56 & 12 & 9 & 0 & \\ \hline
\end{tabular}
\end{center}
\vspace{0.5\baselineskip}
Thus: $\tn^{44}_{0}(001) = 2$, $\tn^{44}_{1}(001) = 56$, 
$\tn^{44}_{2}(001) = 12$. The first corresponds to two sequences: 
$\{ 01010101 \}$, $\{ 10101010 \}$.

Taking in \eqref{e-33b} the sum over $m$, we have also the distribution 
of the string $(001)$ for all $2^{N}$ binary sequences of cardinality 
$N$:
\begin{equation}\label{e-34}
    \begin{split}
        \tn^{N}_{\ell}(001) 
        &= N \sum_{m=\ell+1}^{N} \sum_{h=\ell}^{(m-\ell)\wedge n} 
           \frac{\ell}{(N-m)(m-h)} \binom{h}{\ell} \binom{m-h}{\ell} 
           \binom{N-m}{h}  \\
        &= \frac{N}{\ell} \sum_{m} \binom{N-m-1}{\ell-1} 
           \sum_{h} \binom{m-h-1}{\ell-1} \binom{N-m-\ell}{h-\ell} .
    \end{split}
\end{equation}

Turning to the second expression in \eqref{e-27}, we can see in 
Figure \ref{f-2} that to each tableau of dimension $\ell$ obtained by 
deleting the first column, there corresponds the bottom $h-\ell$ positions 
of $h$ which are responsible for the strings $(101)$, then:

\begin{untitled}
    The number of sequences in which the string $(01)$ appears $h$ 
    times and the string $(101)$ appears $\ell$ times is given by:
    \begin{equation}\label{e-35}
        \tn^{mn}_{h \ell}(01;101) 
        = \tn^{mn}_{h,h-\ell}(01;101) 
        = \frac{N}{h} c^{m}_{h(h-\ell)} \mn^{h}_{n} 
        = \frac{N}{n} c^{m}_{h(h-\ell)} \binom{n}{h}    .
    \end{equation}
\end{untitled}
\noindent Moreover:
\begin{subequations}
\begin{equation}\label{e-36a}
    \tn^{mn}_{\ell}(101) 
    = \frac{N}{n} \sum_{h} c^{m}_{h(h-\ell)} \binom{n}{h}   .
\end{equation}
Here, the index $\ell$ runs from $0$ to $m\wedge(n+\delta_{mn}-1)$ 
because the extreme cases have the form (for $m>n$, $m<n$ and $m=n$):
\[
    \{ 000 \cdots 0101 \cdots 1010 \cdots 000 \} , \quad 
    \{ 111 \cdots 10101 \cdots 101 \cdots 111 \} \quad\text{and}\quad 
    \{ 10101 \cdots 1010 \} .
\]
In the sum, index $h$ goes from $h=\ell$ to 
$h=\left[ \frac{m+\ell}{2} \right]$ (because $h=\ell=m$ gives 
$c^{m}_{m0}=1$ and $m-h=h-\ell$ for the high index), then:

\begin{untitled}
    The number of sequences of $m$ $0$'s and $n$ $1$'s which have 
    $\ell$ strings $(101)$ is:
    \begin{equation}\label{e-36b}
        \begin{split}
            \tn^{mn}_{\ell}(101)
            &= \frac{N}{n} 
               \sum_{h=\ell}^{\left[ \frac{m+\ell}{2} \right] \wedge n} 
               \frac{h-\ell}{m-h} \binom{h}{\ell} \binom{m-h}{h-\ell} 
               \binom{n}{h}    \\
            &= \frac{N}{n} \binom{n}{\ell} 
               \sum_{h=\ell}^{\left[ \frac{m+\ell}{2} \right] \wedge n} 
               \binom{m-h-1}{h-\ell-1} \binom{n-\ell}{h-\ell}   .
        \end{split}
    \end{equation}
\end{untitled}
\end{subequations}

There remains only the last expression in \eqref{e-27}; however we 
must first see how the numbers $\tn(000)$ are related to string 
$(0001)$ and so we must pospone its derivation until after the case 
$|U|=4$. To justify this delay and find a way to obtain $\tn(000)$ 
we will see that, though more complicated than \eqref{e-20}, formula 
\eqref{e-21} for $\tn(00)$ can be found also in terms of strings in 
which $|U|=3$. Let $g$ be a fixed number of strings $(00)$, it is 
easy to see (Figure \ref{f-2}) that $m-h=g$; i.e., there are as many $(00)$ as 
zeros in the subtableau of length $\ell$. After deleting the first 
column of length $h=m-g$ there appears a number $c^{m}_{m-g,\ell}$ of 
such subtableaux, in total an amount of $g c^{m}_{m-g,\ell}$ zeros, 
then:
\[
    \begin{split}
        \tn^{mn}_{g}(00) 
        &= \sum_{\ell=1}^{m-g} \frac{N}{m-g} g c^{m}_{m-g,\ell} 
           \mn^{m-g}_{n} 
         = \sum_{\ell=1}^{m-g} \frac{N}{m-g} g \frac{\ell}{g} 
           \binom{m-g}{\ell} \binom{g}{\ell} \frac{m-g}{n}
           \binom{n}{m-g}    \\
        &= \frac{N}{n} \binom{n}{m-g} \sum_{\ell=0}^{m-g} \ell 
           \binom{m-g}{\ell} \binom{g}{\ell} 
         = \frac{N}{mn} (m-g) \binom{m}{g} \binom{n}{m-g}.
    \end{split}
\]
The last derivation is due to:
\[
    \sum_{0}^{m-g} \ell \binom{g}{\ell} \binom{m-g}{\ell} 
    = \sum_{0}^{g} \ell \binom{g}{\ell} \binom{m-g}{\ell} 
    = \sum_{0}^{(m-g)\wedge g} \ell\binom{g}{\ell}\binom{m-g}{\ell} , 
\]
allowing to use formula \eqref{e-7}.
\section{The Numbers $\tn^{mn}_{\ell}%
         (00\protect\stackrel{s+2}{\cdots}01)$.\\%
         Generalized Kaplansky Lemma}\label{s-TNTGKL}

Define:
\begin{subequations}\label{e-37}
    \begin{align}
        \begin{split}
            c'{}^{i}_{jk} 
            &= \sum_{f=k}^{j \wedge (i-j-k)} \frac{k}{i-j-f} 
               \binom{j}{f} \binom{f}{k} \binom{i-j-f}{k} \\
            &= \binom{j}{k} \sum_{f=k}^{j \wedge (i-j-k)} 
               \binom{j-k}{f-k} \binom{i-j-f-1}{k-1} \quad 
               \text{for $k \neq 0$} ,
        \end{split} \\
\intertext{and:}
            c'{}^{i}_{j0} 
            &= \begin{cases}
                   \binom{j}{i-j} & \text{if $i \leq 2j$} , \\
                   0              & \text{if $i >    2j$} .
               \end{cases}
    \end{align}
\end{subequations}
then,

\begin{lemmaonetwo}
After deleting the first and second columns in all Young-De Moivre tableaux 
of weight $m$ and dimension $h$ one obtains $\mn$-partitions of dimensions 
$\ell = 0,1,2,\ldots,h \wedge (m-2h)$ which number is $c'{}^{m}_{h\ell}$.
\end{lemmaonetwo}

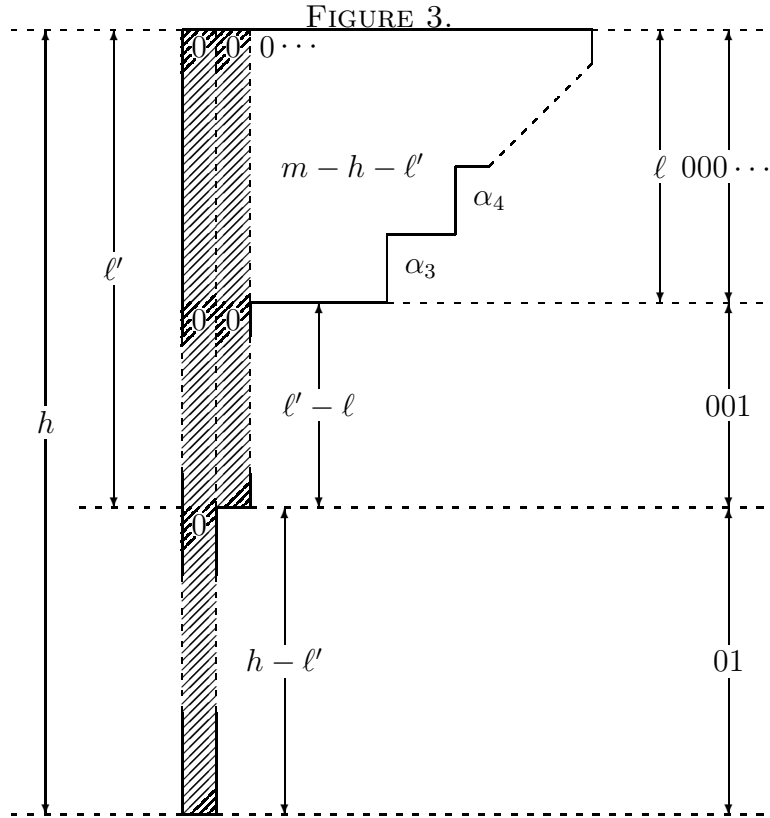
\begin{figure}[H]\caption{}\label{f-3}
\begin{center}
\setlength{\unitlength}{2.5ex}
\begin{picture}(22,23)
\dashline{0.2}(0,0)(5,0)\drawline(5,0)(6,0)\dashline{0.2}(6,0)(22,0)
\dashline{0.2}(2,9)(5,9)\drawline(6,9)(7,9)\dashline{0.2}(7,9)(22,9)
\drawline(7,15)(11,15)\dashline{0.2}(11,15)(22,15)
\drawline(11,17)(13,17)
\drawline(13,19)(14,19)
\dashline{0.2}(0,23)(5,23)\drawline(5,23)(17,23)\dashline{0.2}(17,23)(22,23)
\drawline(5,0)(5, 3)\dashline{0.2}(5, 3)(5, 7)
\drawline(5,7)(5,10)\dashline{0.2}(5,10)(5,14)
\drawline(5,14)(5,23)
\drawline(6,0)(6,3)\dashline{0.2}(6,3)(6, 7)
\drawline(6,7)(6,9)\dashline{0.2}(6,9)(6,23)
\drawline(7, 9)(7,10)\dashline{0.2}(7,10)(7,14)
\drawline(7,14)(7,15)\dashline{0.2}(7,15)(7,23)
\drawline(11,15)(11,17)
\drawline(13,17)(13,19)
\drawline(17,22)(17,23)
\dashline{0.2}(14,19)(17,22)
\put(1,11.5){\makebox(0,0){$h$}}
\put(1,12){\vector(0,1){11}}\put(1,11){\vector(0,-1){11}}
\put(3,16){\makebox(0,0){$\ell'$}}
\put(3,16.5){\vector(0,1){6.5}}\put(3,15.5){\vector(0,-1){6.5}}
\put(8,4.5){\makebox(0,0){$h - \ell'$}}
\put(8,5){\vector(0,1){4}}\put(8,4){\vector(0,-1){4}}
\put(9,12){\makebox(0,0){$\ell' - \ell$}}
\put(9,12.5){\vector(0,1){2.5}}\put(9,11.5){\vector(0,-1){2.5}}
\put(12,16){\makebox(0,0){$\alpha_{3}$}}
\put(14,18){\makebox(0,0){$\alpha_{4}$}}
\put(19,19){\makebox(0,0){$\ell$}}
\put(19,19.5){\vector(0,1){3.5}}\put(19,18.5){\vector(0,-1){3.5}}
\put(21,19){\makebox(0,0){$000\cdots$}}
\put(21,19.5){\vector(0,1){3.5}}\put(21,18.5){\vector(0,-1){3.5}}
\put(21,12){\makebox(0,0){$001$}}
\put(21,12.5){\vector(0,1){2.5}}\put(21,11.5){\vector(0,-1){2.5}}
\put(21,4.5){\makebox(0,0){$01$}}
\put(21,5){\vector(0,1){4}}\put(21,4){\vector(0,-1){4}}
\put(5.5, 8.5){\makebox(0,0){$0$}}
\put(5.5,14.5){\makebox(0,0){$0$}}
\put(6.5,14.5){\makebox(0,0){$0$}}
\put(5.5,22.5){\makebox(0,0){$0$}}
\put(6.5,22.5){\makebox(0,0){$0$}}
\put(7.5,22.5){\makebox(0,0){$0$}}
\put(8.5,22.5){\makebox(0,0){$\cdots$}}
\put(10,19){\makebox(0,0){$m - h - \ell'$}}
%
\multiput(5, 0.0)(0,0.25){31}{\drawline(0,0)(1,1)}
\multiput(5, 8.5)(0,0.25){21}{\drawline(0,0)(1,1)}
\multiput(5,14.5)(0,0.25){29}{\drawline(0,0)(1,1)}
\drawline(5.25,0)(6,0.75)
\drawline(5.50,0)(6,0.50)
\drawline(5.75,0)(6,0.25)
\drawline(5.000,22.750)(5.250,23.000)    \drawline(6.000,22.750)(6.250,23.000)
\drawline(5.000,22.500)(5.500,23.000)    \drawline(6.000,22.500)(6.500,23.000)
\drawline(5.000,22.250)(5.250,22.500)    \drawline(6.000,22.250)(6.250,22.500)
\drawline(5.000,22.000)(5.250,22.250)    \drawline(6.000,22.000)(6.250,22.250)
\drawline(5.000,21.750)(5.375,22.125)    \drawline(6.000,21.750)(6.375,22.125)
\drawline(5.750,22.750)(6.000,23.000)    \drawline(6.750,22.750)(7.000,23.000)
\drawline(5.750,22.500)(6.000,22.750)    \drawline(6.750,22.500)(7.000,22.750)
\drawline(5.625,22.875)(5.750,23.000)    \drawline(6.625,22.875)(6.750,23.000)
\drawline(5.000,14.750)(5.250,15.000)    \drawline(6.000,14.750)(6.250,15.000)
\drawline(5.000,14.500)(5.500,15.000)    \drawline(6.000,14.500)(6.500,15.000)
\drawline(5.000,14.250)(5.250,14.500)    \drawline(6.000,14.250)(6.250,14.500)
\drawline(5.000,14.000)(5.250,14.250)    \drawline(6.000,14.000)(6.250,14.250)
\drawline(5.000,13.750)(5.375,14.125)    \drawline(6.000,13.750)(6.375,14.125)
\drawline(5.750,14.750)(6.000,15.000)    \drawline(6.750,14.750)(7.000,15.000)
\drawline(5.750,14.500)(6.000,14.750)    \drawline(6.750,14.500)(7.000,14.750)
\drawline(5.625,14.875)(6.000,15.250)    \drawline(6.625,14.875)(7.000,15.250)
\drawline(5.000, 8.750)(5.250, 9.000)
\drawline(5.000, 8.500)(5.500, 9.000)
\drawline(5.000, 8.250)(5.250, 8.500)
\drawline(5.000, 8.000)(5.250, 8.250)
\drawline(5.000, 7.750)(5.375, 8.125)
\drawline(5.750, 8.750)(6.000, 9.000)
\drawline(5.750, 8.500)(6.000, 8.750)
\drawline(5.625, 8.875)(6.000, 9.250)
\multiput(6, 9.0)(0,0.25){19}{\drawline(0,0)(1,1)}
\multiput(6,14.5)(0,0.25){29}{\drawline(0,0)(1,1)}
\drawline(6.25,9)(7,9.75)\drawline(6.50,9)(7,9.50)\drawline(6.75,9)(7,9.25)
\drawline(6,22.5)(6.5,23)\drawline(6,22.75)(6.25,23)
%
\end{picture}
\end{center}
\end{figure}

The proof is like I${}^{(1)}$, (see Figure \ref{f-3}).
\noindent The desired number is, (the sum without indexes is extended of all 
$\mn$-partitions):
\[
    \mathcal{N} 
    = \sum_{\ell'} \sum \frac{h!}{(h-\ell')!\,(\ell'-\ell)!} 
      \cdot \frac{1}{\alpha_{3}!\,\alpha_{4}!\,\cdots} 
    = \ell!\, \sum_{\ell'} \binom{h}{\ell'} \binom{\ell'}{\ell} 
      \sum \frac{1}{\alpha_{3}!\,\alpha_{4}!\,\cdots} .
\]
on the other hand:
\[
    \sum \frac{\ell!}{\alpha_{3}!\,\alpha_{4}!\,\cdots} 
    = \ell!\, \sum \frac{1}{\alpha_{3}!\,\alpha_{4}!\,\cdots} 
    = \mn^{\ell}_{m-h-\ell'}  .
\]
Then:
\[
    \mathcal{N} 
    = \sum_{\ell'=\ell}^{h\wedge(m-h-\ell)} 
      \binom{h}{\ell'} \binom{\ell'}{\ell} \mn^{\ell}_{m-h-\ell'}
    = c'{}^{m}_{h\ell} \qquad (\ell \neq 0)  ;
\]
the upper limit is due to the number of elements $m-h-\ell'$ in the 
subtableau which must be at least equal to $\ell$. In the case $\ell=0$ 
the symbol $\mn^{0}_{m-h-\ell'}$ is meaningless. In this case however 
we observe that the condition $m-h-\ell'=0$ gives 
$\mathcal{N} = c'{}^{m}_{h0} = \binom{h}{m-h}$.

\begin{subequations}
\begin{untitled}
    The number of sequences with $m$ zeros and $n$ ones in which the 
    string $(01)$ appears $h$ times, string $(001)$ appears $\ell'$ times 
    and $(0001)$ $\ell$ times is:
    \begin{equation}\label{e-38a}
        \begin{split}
            \tn^{mn}_{h\ell'\ell}(01;001;0001)
            &= \frac{N}{n} \binom{h}{\ell} \binom{h-\ell}{\ell'-\ell} 
               \binom{m-h-\ell'-1}{\ell-1} \binom{n}{h} \\
            &= \frac{N}{n} \binom{n}{\ell} \binom{n-\ell}{\ell'-\ell}
               \binom{n-\ell'}{h-\ell'} \binom{m-h-\ell'-1}{\ell-1} , 
        \end{split}
    \end{equation}
\end{untitled}
\noindent And in case $m \leq 2h$:
\begin{equation}\label{e-38b}
    \tn^{mn}_{h,m-h,0}(01;001;0001) 
        = \frac{N}{n} \binom{n}{m-h} \binom{n-m-h}{2h-m}  ,
\end{equation}
\end{subequations}
which, of course, coincides with 
$\tn^{mn}_{h,m-h}(01;001) = \frac{N}{n} c^{m}_{h,m-h} \binom{n}{h}$, 
because $c^{i}_{j,i-j} = \binom{j}{i-j}$. Though binomial 
coefficients make elementary computations, it is always convenient to 
give fundamental formulas in terms of $c$'s because its matrices can be 
tabulated independently of particular cases. We have for $\tn(0001)$'s:
\begin{subequations}\label{e-39}
    \begin{align}
        \begin{split}
            \tn^{mn}_{\ell}(0001)
            &= \frac{N}{n} \sum c'{}^{m}_{h\ell} \mn^{h}_{n} \\
            &= \frac{N}{n} \binom{n}{\ell} 
               \sum_{h=\ell}^{(m-2\ell)\wedge n} \sum_{\ell'} 
               \binom{n-\ell'}{h-\ell'} \binom{n-\ell}{\ell'-\ell} 
               \binom{m-h-\ell'-1}{\ell-1} \qquad (\ell \neq 0),
        \end{split} \\
        \tn^{mn}_{0}(0001) 
        &= \frac{N}{n} 
           \sum_{h=\left[ \frac{m}{2} \right]}^{n \wedge (m-n)} 
           \binom{n}{m-h} \binom{n-m+h}{2h-m}  .
    \end{align}
\end{subequations}
For the general case we define:
\begin{subequations}\label{e-40}
    \begin{align}
        \begin{split}
            c^{(s)i}_{jk} 
            &= \sum \frac{k}{i-j-f_{1}-f_{2}-\cdots-f_{s}}
               \binom{j}{f_{1}} \binom{f_{1}}{f_{2}} \cdots 
               \binom{f_{s}}{k} \\
            &\quad \cdot \binom{i-j-f_{1}-f_{2}-\cdots-f_{s}}{k} \\
            &= \binom{j}{k} 
               \sum_{j\geq f_{1}\geq f_{2}\geq\cdots\geq f_{s}\geq k} 
               \binom{j-f_{2}}{f_{1}-f_{2}} 
               \binom{j-f_{3}}{f_{2}-f_{3}} \cdots 
               \binom{j-f_{s}}{f_{s-1}-f_{s}} 
               \binom{j-k}{f_{s}-k} \\
            &\quad \cdot \binom{i-j-f_{1}-f_{2}-\cdots-f_{s}-1}{k-1} 
                   \qquad (k \neq 0) ,
        \end{split} \\
        c^{(s)i}_{j0}
        &= \begin{cases}
               \sum_{j\ge f_{1}\ge f_{2}\ge\cdots\ge f_{s}} 
               \binom{j-f_{2}}{f_{1}-f_{2}} \binom{j-f_{3}}{f_{2}-f_{3}} 
               \cdots \binom{j-f_{s}}{f_{s-1}-f_{s}} 
               \binom{j}{f_{s}} & \text{if $i \leq (s+1)j$} , \\ 
               0                & \text{if $i >    (s+1)j$} ;
           \end{cases}
    \end{align}
\end{subequations}
then,

\begin{lemmaone}
After deleting the $s+1$ columns in all Young-De Moivre 
tableaux of weight $m$ and dimension $h$ one obtain $\mn$-partitions of 
dimensions $\ell = 0,1,2,\ldots,h \wedge (m-(s+1)h)$ which number is 
$c^{(s)m}_{h\ell}$. \textup{(This generalizes I${}^{(1)}$ and I${}^{(2)}$.)}
\end{lemmaone}

\noindent So, we have:
\[
    \tn^{mn}_{h \ell' \ell'' \cdots \ell^{(s)}}
        (01;001;\cdots;00\stackrel{(s+2)}{\cdots}01)
    = \frac{N}{h} \overline{c}{}^{(s)m}_{h\ell} \mn^{h}_{n} ,
\]
$\overline{c}{}^{(s)m}_{h\ell}$ being the $c^{(s)m}_{h\ell}$ without 
summation on $\ell',\ell'',\ldots,\ell^{(s)}$. And:
\begin{equation}\label{e-41}
    \begin{split}
        \tn^{mn}_{\ell}(00 \cdots 01) 
        &= \frac{N}{n} \binom{n}{\ell} 
           \sum_{h=\ell}^{(m-(s+1)\ell)\wedge n} 
           \sum_{\ell'\geq \ell''\geq \cdots\geq \ell^{(s)}\geq h}
           \binom{h-\ell''}{\ell'-\ell''} 
           \binom{h-\ell'''}{\ell''-\ell'''} \cdots 
           \binom{h-\ell^{(s)}}{\ell^{(s-1)}-\ell^{(s)}} \\
        &\quad \cdot \binom{n-\ell}{h-\ell} 
               \binom{m-h-\ell'-\ell''-\cdots-\ell^{(s)}-1}{\ell-1} . 
    \end{split}
\end{equation}
Formula \eqref{e-41} generalizes Kaplansky's lemma (for a redefinition 
of Kaplansky's lemma see Comptet \cite[p.~35]{Cl1970}). In particular, 
taking $h = \ell' = \ell'' = \cdots = \ell^{(s)} = n$, we have:
\[
    g_{n}(N,p) = \frac{N}{n} \binom{m-n-sn-1}{n-1} 
               = \frac{N}{N-pn} \binom{N-pn}{n} , \tag{Kaplansky}
\]
being $p=s+2$ the least length of blocks which separate the ones.
\section{The Numbers $\tn^{mn}_{g}%
         (00\protect\stackrel{s+2}{\cdots}00)$.\\%
         Fibonacci $\gfn(N;r,h)$ Numbers}\label{s-TNTFNF}

Our final problem will consist of looking at the weights $g$ of the 
subtableaux, instead of the dimension $\ell$. We will call 
$c^{(s)mg}_{h}$ the number of subtableaux of weight $g$, which remain 
after deleting $s+1$ succesive columns, starting from the first, in 
all the Young-De Moivre tableaux of weight $m~(>g)$ and dimension $h$.

For $s=0$, in Figure \ref{f-2} we can see that the weight of the remaining 
subtableaux is $m-h$, then (\textbf{II}$\boldsymbol{{}^{(1)}}$):
\begin{equation}\label{e-42}
    \begin{split}
        c^{mg}_{h} 
        &= \sum_{\ell=0}^{h} c^{m}_{m-g,\ell} 
         = \sum_{0}^{m-g} \frac{\ell}{g} \binom{m-g}{\ell} \binom{g}{\ell} 
         = \frac{m-g}{m} \binom{m}{g} \\
        &= \frac{m-g}{m} \binom{m}{m-g} 
         = \binom{m-1}{g} 
         = c^{mg}
    \end{split}
\end{equation}
so $c^{mg}_{h} = c^{mg}$ does not depend on $h$, i.e.: Counting shall 
be made on all tableaux of weight $m$. Though indirectly, it is thanks 
to numbers $c^{mg}_{h}$ that $\tn(00)$ was found and we may now write 
formula \eqref{e-21} as $\frac{N}{m-g} c^{mg} \mn^{m-g}_{n}$.

For $s=1$, it is easy to see that the problem of counting strings 
$(000)$ is given also by the formula \eqref{e-39} but collecting the 
terms with the same $g$, i.e., from \eqref{e-37} we obtain 
$c'{}^{ig}_{j}$ taking in the sum only the term $f=i-j-g$:
\[
    c'{}^{i}_{jk}(g) = \frac{k}{g} \binom{j}{i-j-g} \binom{i-j-g}{k} 
                       \binom{g}{k} , \qquad
    c'{}^{ig}_{j} = \sum_{k} c'{}^{i}_{jk}(g) .
\]
We have (see Figure \ref{f-3}):
\begin{subequations}\label{e-43}
    \begin{align}
        \begin{split}
            c'{}^{mg}_{h} 
            &= \sum_{\ell} \sum_{m-h-\ell'=g} \binom{h}{\ell'} 
               \binom{\ell'}{\ell} \mn^{\ell}_{m-h-\ell'} 
             = \sum_{\ell=0}^{h \wedge g} \binom{h}{\ell} 
                \binom{h-\ell}{m-g-h-\ell} \mn^{\ell}_{g} \\
            &= \frac{1}{g}  \sum_{\ell=0}^{h \wedge g \wedge (m-h-g)} 
               \ell \binom{h}{\ell} \binom{g}{\ell} 
               \binom{h-\ell}{m-g-h-\ell} , \qquad (g \neq 0) ;
        \end{split} \\
        c'{}^{m0}_{h} 
        &= c'{}^{m}_{h0} 
         = \begin{cases}
               \binom{h}{m-h} &\text{if $m \leq 2h$} , \\
               0              &\text{if $m >    2h$} . 
           \end{cases}
    \end{align}
\end{subequations}
Then we have:

\begin{lemmatwotwo}
After deleting the first and second column in all Young-De Moivre tableaux of 
weight $m$ and dimension $h$ one obtains $\mn$-partitions of weights 
$g = 0 \vee (m-2h), 1 \vee (m-2h), \ldots, m-h-1$ which number is 
$c'{}^{mg}_{h}$.
\end{lemmatwotwo}

In order to clarify ambiguity, let us make the process for $m=7$, 
$h=3$:
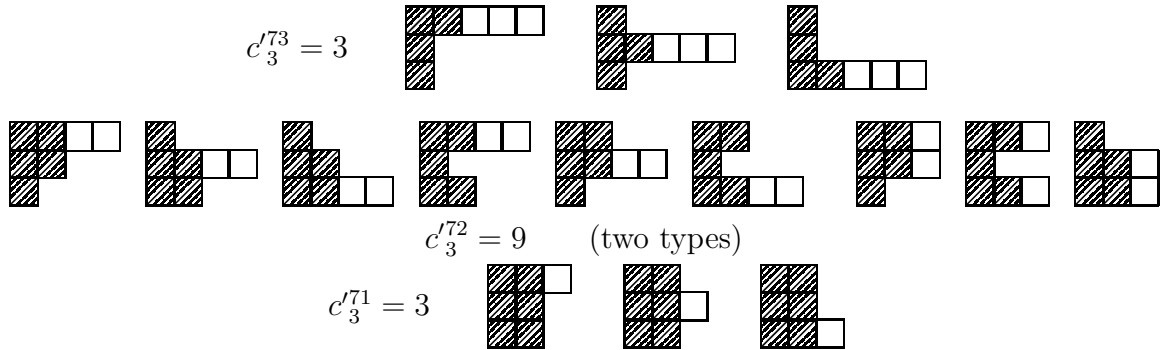
\setlength{\unitlength}{2ex}
\savebox{\shbx}(1,1){\thinlines%
                     \drawline(0,0.75)(0.25,1)%
                     \drawline(0,0.50)(0.50,1)%
                     \drawline(0,0.25)(0.75,1)%
		     \drawline(0,0)(1,1)%
 		     \drawline(0.25,0)(1,0.75)%
		     \drawline(0.50,0)(1,0.50)%
                     \drawline(0.75,0)(1,0.25)}%
\begin{center}
\begin{picture}(25,3)(-6,0)
%
\put(-6,1){\makebox(4,1){$c'{}^{73}_{3} = 3$}}
\multiput( 0,2)(1,0){5}{\framebox(1,1){}}\multiput( 0,2)(1,0){2}{\usebox{\shbx}}
\multiput( 0,1)(1,0){1}{\framebox(1,1){}}\multiput( 0,1)(1,0){1}{\usebox{\shbx}}
\multiput( 0,0)(1,0){1}{\framebox(1,1){}}\multiput( 0,0)(1,0){1}{\usebox{\shbx}}
\multiput( 7,2)(1,0){1}{\framebox(1,1){}}\multiput( 7,2)(1,0){1}{\usebox{\shbx}}
\multiput( 7,1)(1,0){5}{\framebox(1,1){}}\multiput( 7,1)(1,0){2}{\usebox{\shbx}}
\multiput( 7,0)(1,0){1}{\framebox(1,1){}}\multiput( 7,0)(1,0){1}{\usebox{\shbx}}
\multiput(14,2)(1,0){1}{\framebox(1,1){}}\multiput(14,2)(1,0){1}{\usebox{\shbx}}
\multiput(14,1)(1,0){1}{\framebox(1,1){}}\multiput(14,1)(1,0){1}{\usebox{\shbx}}
\multiput(14,0)(1,0){5}{\framebox(1,1){}}\multiput(14,0)(1,0){2}{\usebox{\shbx}}
\end{picture}
\end{center}
\vspace{0.5\baselineskip}
\begin{center}
\begin{picture}(42,3)
%
\multiput( 0,2)(1,0){4}{\framebox(1,1){}}\multiput( 0,2)(1,0){2}{\usebox{\shbx}}
\multiput( 0,1)(1,0){2}{\framebox(1,1){}}\multiput( 0,1)(1,0){2}{\usebox{\shbx}}
\multiput( 0,0)(1,0){1}{\framebox(1,1){}}\multiput( 0,0)(1,0){1}{\usebox{\shbx}}
\multiput( 5,2)(1,0){1}{\framebox(1,1){}}\multiput( 5,2)(1,0){1}{\usebox{\shbx}}
\multiput( 5,1)(1,0){4}{\framebox(1,1){}}\multiput( 5,1)(1,0){2}{\usebox{\shbx}}
\multiput( 5,0)(1,0){2}{\framebox(1,1){}}\multiput( 5,0)(1,0){2}{\usebox{\shbx}}
\multiput(10,2)(1,0){1}{\framebox(1,1){}}\multiput(10,2)(1,0){1}{\usebox{\shbx}}
\multiput(10,1)(1,0){2}{\framebox(1,1){}}\multiput(10,1)(1,0){2}{\usebox{\shbx}}
\multiput(10,0)(1,0){4}{\framebox(1,1){}}\multiput(10,0)(1,0){2}{\usebox{\shbx}}
\multiput(15,2)(1,0){4}{\framebox(1,1){}}\multiput(15,2)(1,0){2}{\usebox{\shbx}}
\multiput(15,1)(1,0){1}{\framebox(1,1){}}\multiput(15,1)(1,0){1}{\usebox{\shbx}}
\multiput(15,0)(1,0){2}{\framebox(1,1){}}\multiput(15,0)(1,0){2}{\usebox{\shbx}}
\multiput(20,2)(1,0){2}{\framebox(1,1){}}\multiput(20,2)(1,0){2}{\usebox{\shbx}}
\multiput(20,1)(1,0){4}{\framebox(1,1){}}\multiput(20,1)(1,0){2}{\usebox{\shbx}}
\multiput(20,0)(1,0){1}{\framebox(1,1){}}\multiput(20,0)(1,0){1}{\usebox{\shbx}}
\multiput(25,2)(1,0){2}{\framebox(1,1){}}\multiput(25,2)(1,0){2}{\usebox{\shbx}}
\multiput(25,1)(1,0){1}{\framebox(1,1){}}\multiput(25,1)(1,0){1}{\usebox{\shbx}}
\multiput(25,0)(1,0){4}{\framebox(1,1){}}\multiput(25,0)(1,0){2}{\usebox{\shbx}}
\multiput(31,2)(1,0){3}{\framebox(1,1){}}\multiput(31,2)(1,0){2}{\usebox{\shbx}}
\multiput(31,1)(1,0){3}{\framebox(1,1){}}\multiput(31,1)(1,0){2}{\usebox{\shbx}}
\multiput(31,0)(1,0){1}{\framebox(1,1){}}\multiput(31,0)(1,0){1}{\usebox{\shbx}}
\multiput(35,2)(1,0){3}{\framebox(1,1){}}\multiput(35,2)(1,0){2}{\usebox{\shbx}}
\multiput(35,1)(1,0){1}{\framebox(1,1){}}\multiput(35,1)(1,0){1}{\usebox{\shbx}}
\multiput(35,0)(1,0){3}{\framebox(1,1){}}\multiput(35,0)(1,0){2}{\usebox{\shbx}}
\multiput(39,2)(1,0){1}{\framebox(1,1){}}\multiput(39,2)(1,0){1}{\usebox{\shbx}}
\multiput(39,1)(1,0){3}{\framebox(1,1){}}\multiput(39,1)(1,0){2}{\usebox{\shbx}}
\multiput(39,0)(1,0){3}{\framebox(1,1){}}\multiput(39,0)(1,0){2}{\usebox{\shbx}}
\end{picture}
\end{center}
\[
    c'{}^{72}_{3} = 9 \qquad \makebox{(two types)}
\]
\begin{center}
\begin{picture}(19,3)(-6,0)
%
\put(-6,1){\makebox(4,1){$c'{}^{71}_{3} = 3$}}
\multiput( 0,2)(1,0){3}{\framebox(1,1){}}\multiput( 0,2)(1,0){2}{\usebox{\shbx}}
\multiput( 0,1)(1,0){2}{\framebox(1,1){}}\multiput( 0,1)(1,0){2}{\usebox{\shbx}}
\multiput( 0,0)(1,0){2}{\framebox(1,1){}}\multiput( 0,0)(1,0){2}{\usebox{\shbx}}
\multiput( 5,2)(1,0){2}{\framebox(1,1){}}\multiput( 5,2)(1,0){2}{\usebox{\shbx}}
\multiput( 5,1)(1,0){3}{\framebox(1,1){}}\multiput( 5,1)(1,0){2}{\usebox{\shbx}}
\multiput( 5,0)(1,0){2}{\framebox(1,1){}}\multiput( 5,0)(1,0){2}{\usebox{\shbx}}
\multiput(10,2)(1,0){2}{\framebox(1,1){}}\multiput(10,2)(1,0){2}{\usebox{\shbx}}
\multiput(10,1)(1,0){2}{\framebox(1,1){}}\multiput(10,1)(1,0){2}{\usebox{\shbx}}
\multiput(10,0)(1,0){3}{\framebox(1,1){}}\multiput(10,0)(1,0){2}{\usebox{\shbx}}
\end{picture}
\end{center}
In $c'{}^{72}_{3}$ tableaux, first and fourth give the same position 
for the subtableau of weight $2$, but since the second column is 
different that will give different matchings; for example ($n=9$),
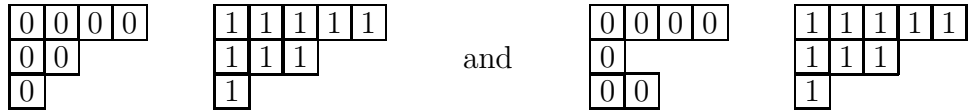
\begin{center}
\setlength{\unitlength}{2.5ex}
\begin{picture}(28,3)
%
\multiput( 0,2)(1,0){4}{\framebox(1,1){0}}
\multiput( 0,1)(1,0){2}{\framebox(1,1){0}}
\multiput( 0,0)(1,0){1}{\framebox(1,1){0}}
\multiput( 6,2)(1,0){5}{\framebox(1,1){1}}
\multiput( 6,1)(1,0){3}{\framebox(1,1){1}}
\multiput( 6,0)(1,0){1}{\framebox(1,1){1}}
\put(13,1){\makebox(2,1){and}}
\multiput(17,2)(1,0){4}{\framebox(1,1){0}}
\multiput(17,1)(1,0){1}{\framebox(1,1){0}}
\multiput(17,0)(1,0){2}{\framebox(1,1){0}}
\multiput(23,2)(1,0){5}{\framebox(1,1){1}}
\multiput(23,1)(1,0){3}{\framebox(1,1){1}}
\multiput(23,0)(1,0){1}{\framebox(1,1){1}}
\end{picture}
\end{center}
give different sequences,
\[
    \{ 0000111110011101 \} \qquad \text{and} \qquad 
    \{ 0000111110111001 \} .
\]

Now, at last, the coefficients $c'{}^{mg}_{h}$ provide for the numbers 
$\tn(000)$ (and so for $\tn(111)$). In fact, there are as many strings 
$(000)$ as zeros ($g$) in the subtableaux of weight $g$, so that:
\begin{equation}\label{e-44}
    \tn^{mn}_{g}(000) 
    = \sum_{h} \frac{N}{h} c'{}^{mg}_{h} \mn^{h}_{n} 
    = \frac{N}{n} \sum_{h} c'{}^{mg}_{h} \binom{n}{h} .
\end{equation}
Since the matrices $c'{}^{mg}_{h}$ ($g$ fixed) can be tabulated very 
easily, as in Appendix A, it is preferable to work with them; 
nonetheless if we write explicitly all binomial coefficients, we have:
\begin{subequations}\label{e-45}
    \begin{align}
        \tn^{mn}_{g}(000)
        &= \frac{N}{n} \sum_{h} \sum_{\ell} \ell \binom{n}{\ell} 
           \binom{n}{h} \binom{n-\ell}{m-h-g-\ell} 
           \binom{n-m+h+g}{2h+g-m} , \quad(g \neq 0)  \\
        \tn^{mn}_{0}(000) 
        &= \frac{N}{n} \sum_{h} \binom{h}{m-h} \binom{n}{h} 
         = \sum_{h} \binom{n}{m-h} \binom{n-m+h}{2h-m} .
    \end{align}
\end{subequations}

Let us find a complete example with $m=5$, $n=3$, ($N=8$) which can be 
easily checked by hand on the $56$ sequences. Appendix A provide for 
all coefficients:
\[
    \begin{split}
        \tn^{53}_{0}(000) &= \frac{8}{3} c'{}^{50}_{3} \binom{3}{3} 
                           = \frac{8}{3} 3 \binom{3}{3} 
                           = 8 ; \\
        \tn^{53}_{1}(000) &= \frac{8}{3} 
                             \left( c'{}^{51}_{2} \binom{3}{2} 
                                  + c'{}^{51}_{3} \binom{3}{3} \right) 
                           = \frac{8}{3} 
                             \left( 2 \binom{3}{2}
                                  + 3 \binom{3}{3} \right) 
                           = 24 ; \\
        \tn^{53}_{2}(000) &= \frac{8}{3} c'{}^{52}_{2} \binom{3}{2} 
                           = 16 ; \\
        \tn^{53}_{3}(000) &= \frac{8}{3} c'{}^{53}_{1} \binom{3}{1} 
                           = 8 ; \\
        \tn^{53}_{0}(111) &= \tn^{35}_{0}(000) 
                           = \frac{8}{5} 
                             \left( c'{}^{30}_{2} \binom{5}{2} 
                                  + c'{}^{30}_{3} \binom{5}{3} \right) 
                           = 48 ; \\
        \tn^{53}_{1}(111) &= \frac{8}{5} c'{}^{31}_{1}  
                           = 8 ; \\
        \tn^{53}_{0}(001) &= 0 ; \\
        \tn^{53}_{1}(001) &= \frac{8}{3} 
                             \left( c^{5}_{11} \binom{3}{1} 
                                  + c^{5}_{21} \binom{3}{2} 
                                  + c^{5}_{31} \binom{3}{3} \right) 
                           = 32 ; \quad \text{etc.} \\
        \tn^{53}_{0}(011) &= \tn^{35}_{0}(100) 
                           = \tn^{35}_{0}(001) 
                           = \frac{8}{5}
                             \left( c^{3}_{10} \binom{5}{1} 
                                  + c^{3}_{20} \binom{5}{2} 
                                  + c^{3}_{30} \binom{5}{3} \right) 
                           = 16 ; \quad \text{etc.}
    \end{split}
\]
The following table gives the full results:
\begin{table}[H]\caption{}\label{t-5}%
\begin{center}%
\begin{tabular}{ l c c c c c c c c }\hline
\tstrut $\tn$          & (111) & (110) & (101) & (100) & (011) & (010) & (001) & (000)\\ \hline
\tstrut $\tn^{53}_{0}$ &    48 &    16 &    24 &  \h 0 &    16 &  \h 8 &  \h 0 &  \h 8\\ \hline
\tstrut $\tn^{53}_{1}$ &  \h 8 &    40 &    24 &    32 &    40 &    32 &    32 &    24\\ \hline
\tstrut $\tn^{53}_{2}$ &  \h 0 &  \h 0 &  \h 8 &    24 &  \h 0 &  \h 0 &    24 &    16\\ \hline
\tstrut $\tn^{53}_{3}$ &  \h 0 &  \h 0 &  \h 0 &  \h 0 &  \h 0 &    16 &  \h 0 &  \h 8\\ \hline
\end{tabular}
\end{center}
\end{table}

In general, for $s \geq 1$, take $i-j-f_{1}-f_{2}-\cdots-f_{s}=g$ and 
define:
\begin{subequations}\label{e-46}
    \begin{align}
        c^{(s)ig}_{j} 
        &= \frac{1}{g} \sum_{k=1}^{g} 
           \sum_{\substack{
                 f_{1}\geq f_{2}\geq\cdots\geq f_{s}\\
                 f_{1} + f_{2} + \cdots + f_{s} = i-j-g }} k
                 \binom{j}{f_{1}} \binom{f_{1}}{f_{2}} \cdots 
                 \binom{f_{s}}{k} \binom{g}{k} , 
                 \quad (j\geq f_{i}\geq k ; g\neq 0), \\
\intertext{and:}
        c^{(s)i0}_{j} 
        &= \sum_{\substack{
                 f_{1}\geq f_{2}\geq\cdots\geq f_{s}\\
                 f_{1} + f_{2} + \cdots + f_{s} = i-j }} 
                 \binom{j}{f_{1}} \binom{f_{1}}{f_{2}} \cdots 
                 \binom{f_{s-1}}{f_{s}} , 
                 \quad (j\geq f_{i}\geq k) .
    \end{align}
\end{subequations}
(From hereon in order to calculate the $c$'s one has to deal with the 
number of partitions $\pn^{s}_{i-j}$, $\pn^{s}_{i-j-g}$, etc.\ which 
complicates the formulas). Then we have:

\begin{lemmatwo}
After deleting the first $s+1$ columns in all Young-De Moivre tableaux of 
weight $m$ and dimension $h$ one obtains $\mn$-partitions of weights 
$g = 0 \vee (m-(s+1)h), 1 \vee (m-(s+1)h), \ldots , m-h-s$ which number 
is $c^{(s)mg}_{h}$.
\textup{(This generalizes II${}^{(1)}$ and II${}^{(2)}$.)}
\end{lemmatwo}

This lemma provides for $U=(00\stackrel{s+2}{\cdots}0)$ or 
$U=(11\stackrel{s+2}{\cdots}1)$:
\begin{subequations}\label{e-47}
    \begin{align}
        \tn^{mn}_{g}(00\stackrel{s+2}{\cdots}0) 
        &= \sum_{h} \frac{N}{h} c^{(s)mg}_{h} \mn^{h}_{n} 
         = \frac{N}{n} \sum_{h} c^{(s)mg}_{h} \binom{n}{h} , \\
\intertext{and:}
        \tn^{mn}_{g}(11\stackrel{s+2}{\cdots}1) 
        &= \frac{N}{m} \sum_{h} c^{(s)ng}_{h} \binom{m}{h} .
    \end{align}
\end{subequations}
So, the generalized Fibonacci numbers are:
\begin{equation}\label{e-48}
    \gfn(N;s,g) 
    = N \sum_{n} \sum_{h} \frac{c^{(s-2)ng}_{h}}{N-n} \binom{N-n}{h} .
\end{equation}
For example, the number of subsets of $X = \{ 1,2,3,\ldots,N \}$ that do 
not contain three consecutive integers modulo $N$ is:
\begin{equation}\label{e-49}
    \gfn(N;3,0) 
    = N \sum_{n}^{N} \sum_{h=\left[ \frac{n}{2} \right]}^{N-n} 
      \frac{c'{}^{n0}_{h}}{N-n} \binom{N-n}{h} 
    = N \sum_{n=1}^{n+\left[ \frac{n}{2} \right] \leq N} 
        \sum_{h=\left[ \frac{n}{2} \right]}^{N-n} 
        \frac{ \binom{h}{n-h} \binom{N-n}{h} }{N-n} .
\end{equation}
\begin{example}
In \eqref{e-49} take $N=5$, we have:
\[
    \begin{split}
        \gfn(5;3,0) 
        &= 5 \left\{ 0 + 0 + \frac{1}{2} \binom{2}{3-2} \binom{2}{2} 
                   + \frac{1}{3} \left[ \binom{1}{1} \binom{3}{1} 
                                      + \binom{2}{0} \binom{3}{2} 
                                 \right] 
                   + \frac{1}{4} \binom{1}{0} \binom{4}{1} \right\} \\
        &= 20 .
    \end{split}
\]
The twenty subsets are:
\[
    \begin{split}
        & \{1,2,4\}, \{1,3,4\}, \{1,3,5\}, \{2,3,5\}, \{2,4,5\}, \\
        & \{1,2\},   \{1,3\},   \{1,4\},   \{1,5\},   \{2,3\}, 
          \{2,4\},   \{2,5\},   \{3,4\},   \{3,5\},              \\
        & \{1\},     \{2\},     \{3\},     \{4\},     \{5\} .
    \end{split}
\]
\end{example}
\section*{Appendix A}\label{a-A}

The recurrence relation \eqref{e-31} gives the matrices $c^{i}_{jk}$, we have, 
($i$ rows, $j$ columns, $k$ fixed):
\vspace{0.5\baselineskip}
\begin{center}$c^{i}_{j0} = ${\footnotesize%
\begin{tabular}{ r | r r r r r r r r r r r r }
       & \K 1 & \K 2 & \K 3 & \K 4 & \K 5 & \K 6 & \K 7 & \K 8 & \K 9 & \h 10 & \h 11 & \h 12 \\
    \hline
     1 & 1 \\
     2 & & 2 \\
     3 & & & 3 \\
     4 & & & & 4 \\
     5 & & & & & 5 \\
     6 & & & & & & 6 \\
     7 & & & & & & & 7 \\
     8 & & & & & & & & 8 \\
     9 & & & & & & & & & 9 \\
    10 & & & & & & & & & & 10 \\
    11 & & & & & & & & & & & 11 \\
    12 & & & & & & & & & & & & 12 
\end{tabular}
    }    
\end{center}
\vspace{0.5\baselineskip}
\begin{center}$c^{i}_{j1} = ${\footnotesize%
\begin{tabular}{ r | r r r r r r r r r r r r }
       & \K 1 & \K 2 & \K 3 & \K 4 & \K 5 & \K 6 & \K 7 & \K 8 & \K 9 & \h 10 & \h 11 & \h 12 \\
    \hline
     1 \\
     2 & 1 \\
     3 & 1 & 2 \\
     4 & 1 & 2 & 3 \\
     5 & 1 & 2 & 3 & 4 \\
     6 & 1 & 2 & 3 & 4 & 5 \\
     7 & 1 & 2 & 3 & 4 & 5 & 6 \\
     8 & 1 & 2 & 3 & 4 & 5 & 6 & 7 \\
     9 & 1 & 2 & 3 & 4 & 5 & 6 & 7 & 8 \\
    10 & 1 & 2 & 3 & 4 & 5 & 6 & 7 & 8 & 9 \\
    11 & 1 & 2 & 3 & 4 & 5 & 6 & 7 & 8 & 9 & 10 \\
    12 & 1 & 2 & 3 & 4 & 5 & 6 & 7 & 8 & 9 & 10 & 11 
\end{tabular}
    }
\end{center}
\vspace{0.5\baselineskip}
\begin{center}$c^{i}_{j2} = ${\footnotesize%
\begin{tabular}{ r | r r r r r r r r r r r r }
       & \K 1 & \K 2 & \K 3 & \K 4 & \K 5 & \K 6 & \K 7 & \K 8 & \K 9 & \h 10 & \h 11 & \h 12 \\
    \hline
     1 & \\
     2 & \\
     3 & \\
     4 & & 1 \\
     5 & & 2 & 3 \\
     6 & & 3 & 6 & 6 \\
     7 & & 4 & 9 & 12 & 10 \\
     8 & & 5 & 12 & 18 & 20 & 15 \\
     9 & & 6 & 15 & 24 & 30 & 30 & 21 \\
    10 & & 7 & 18 & 30 & 40 & 45 & 42 & 28 \\
    11 & & 8 & 21 & 36 & 50 & 60 & 63 & 56 & 36 \\
    12 & & 9 & 24 & 42 & 60 & 75 & 84 & 84 & 72 & 45 
\end{tabular}
    }
\end{center}
\vspace{0.5\baselineskip}
\begin{center}$c^{i}_{j3} = ${\footnotesize%
\begin{tabular}{ r | r r r r r r r r r r r r }
       & \K 1 & \K 2 & \K 3 & \K 4 & \K 5 & \K 6 & \K 7 & \K 8 & \K 9 & \h 10 & \h 11 & \h 12 \\
    \hline
     1 & & \\
     2 & & \\
     3 & & \\
     4 & & \\
     5 & & \\
     6 & & & 1 \\
     7 & & & 3 & 4 \\
     8 & & & 6 & 12 & 10 \\
     9 & & & 10 & 24 & 30 & 20 \\
    10 & & & 15 & 40 & 60 & 60 & 35 \\
    11 & & & 21 & 60 & 100 & 120 & 105 & 56 \\
    12 & & & 28 & 84 & 150 & 200 & 210 & 168 & 84 
\end{tabular}
    }
\end{center}
\vspace{0.5\baselineskip}
\begin{center}$c^{i}_{j4} = ${\footnotesize%
\begin{tabular}{ r | r r r r r r r r r r r r }
       & \K 1 & \K 2 & \K 3 & \K 4 & \K 5 & \K 6 & \K 7 & \K 8 & \K 9 & \h 10 & \h 11 & \h 12 \\
    \hline
     1 & & & \\
     2 & & & \\
     3 & & & \\
     4 & & & \\
     5 & & & \\
     6 & & & \\
     7 & & & \\
     8 & & & & 1 \\
     9 & & & & 4 & 5 \\
    10 & & & & 10 & 20 & 15 \\
    11 & & & & 20 & 50 & 60 & 35 \\
    12 & & & & 35 & 100 & 150 & 140 & 70
\end{tabular}
    }
\end{center}
\vspace{0.5\baselineskip}
\begin{center}$c^{i}_{j5} = ${\footnotesize%
\begin{tabular}{ r | r r r r r r r r r r r r }
       & \K 1 & \K 2 & \K 3 & \K 4 & \K 5 & \K 6 & \K 7 & \K 8 & \K 9 & \h 10 & \h 11 & \h 12 \\
    \hline
     1 & & & & \\
     2 & & & & \\
     3 & & & & \\
     4 & & & & \\
     5 & & & & \\
     6 & & & & \\
     7 & & & & \\
     8 & & & & \\
     9 & & & & \\
    10 & & & & & 1 \\
    11 & & & & & 5 & 6 \\
    12 & & & & & 15 & 30 & 21 %
\end{tabular}
    }
\end{center}
\vspace{0.5\baselineskip}

And the finite matrices $c^{i}_{jk}$ ($i$ fixed, $j$ rows, $k$ columns) are:

\vspace{0.5\baselineskip}
\noindent
\begin{minipage}{0.5\textwidth}%
\begin{center}%
{$c^{3}_{jk} = $}
    {\footnotesize%
\begin{tabular}{ r | r r r r r r }%
       & \h 0 & \h 1 & \h 2 & \h 3 & \h 4 & \h 5 \\
    \hline%
     1 & & 1 \\%
     2 & & 2 \\%
     3 & 3 \\%
     4 & \\%
     5 & \\%
     6 & \\%
     7 & \\%
     8 & \\%
     9 & \\%
    10 & %
\end{tabular}%
    }
\end{center}%
\end{minipage}%
\hfill
\begin{minipage}{0.5\textwidth}%
\begin{center}%
{$c^{4}_{jk} = $}
    {\footnotesize%
\begin{tabular}{ r | r r r r r r }%
       & \h 0 & \h 1 & \h 2 & \h 3 & \h 4 & \h 5 \\
    \hline%
     1 & & 1 \\%
     2 & & 2 & 1 \\%
     3 & & 3 \\%
     4 & 4 \\%
     5 & \\%
     6 & \\%
     7 & \\%
     8 & \\%
     9 & \\%
    10 & %
\end{tabular}%
    }
\end{center}%
\end{minipage}%
\vspace{0.5\baselineskip}
\noindent
\begin{minipage}{0.5\textwidth}%
\begin{center}%
{$c^{5}_{jk} = $}
    {\footnotesize%
\begin{tabular}{ r | r r r r r r }
       & \h 0 & \h 1 & \h 2 & \h 3 & \h 4 & \h 5 \\
    \hline
     1 & & 1 \\
     2 & & 2 & 2 \\
     3 & & 3 & 3 \\
     4 & & 4 \\
     5 & 5 \\
     6 & \\
     7 & \\
     8 & \\
     9 & \\
    10 & %
\end{tabular}%
    }
\end{center}%
\end{minipage}%
\hfill
\begin{minipage}{0.5\textwidth}%
\begin{center}%
{$c^{6}_{jk} = $}
    {\footnotesize%
\begin{tabular}{ r | r r r r r r }
       & \h 0 & \h 1 & \h 2 & \h 3 & \h 4 & \h 5 \\
    \hline
     1 & & 1 \\
     2 & & 2 & 3 \\
     3 & & 3 & 5 & 1 \\
     4 & & 4 & 6 \\
     5 & & 5 \\
     6 & 6 \\
     7 & \\
     8 & \\
     9 & \\
    10 & %
\end{tabular}%
    }
\end{center}%
\end{minipage}%
\vspace{0.5\baselineskip}
\noindent
\begin{minipage}{0.5\textwidth}%
\begin{center}%
{$c^{7}_{jk} = $}
    {\footnotesize%
\begin{tabular}{ r | r r r r r r }
       & \h 0 & \h 1 & \h 2 & \h 3 & \h 4 & \h 5 \\
    \hline
     1 & & 1 \\
     2 & & 2 & 4 \\
     3 & & 3 & 9 & 3 \\
     4 & & 4 & 12 & 4 \\
     5 & & 5 & 10 \\
     6 & & 6 \\
     7 & 7 \\
     8 & \\
     9 & \\
    10 & %
\end{tabular}%
    }
\end{center}%
\end{minipage}%
\hfill
\begin{minipage}{0.5\textwidth}%
\begin{center}%
{$c^{8}_{jk} = $}
    {\footnotesize%
\begin{tabular}{ r | r r r r r r }
       & \h 0 & \h 1 & \h 2 & \h 3 & \h 4 & \h 5 \\
    \hline
     1 & & 1 \\
     2 & & 2 & 5 \\
     3 & & 3 & 12 & 6 \\
     4 & & 4 & 18 & 12 & 1 \\
     5 & & 5 & 20 & 10 \\
     6 & & 6 & 15 \\
     7 & & 7 \\
     8 & 8 \\
     9 & \\
    10 & %
\end{tabular}%
    }
\end{center}%
\end{minipage}%
\vspace{0.5\baselineskip}
\noindent%
\begin{minipage}{0.5\textwidth}%
\begin{center}%
{$c^{9}_{jk} = $}
    {\footnotesize%
\begin{tabular}{ r | r r r r r r }
       & \h 0 & \h 1 & \h 2 & \h 3 & \h 4 & \h 5 \\
    \hline
     1 & & 1 \\
     2 & & 2 & 6 \\
     3 & & 3 & 15 & 10 \\
     4 & & 4 & 24 & 24 & 4 \\
     5 & & 5 & 30 & 30 & 5 \\
     6 & & 6 & 30 & 20 \\
     7 & & 7 & 21 \\
     8 & & 8 \\
     9 & 9 \\
    10 & %
\end{tabular}%
    }
\end{center}%
\end{minipage}%
\hfill
\begin{minipage}{0.5\textwidth}%
\begin{center}%
{$c^{10}_{jk} = $}
    {\footnotesize%
\begin{tabular}{ r | r r r r r r }
       & \h 0 & \h 1 & \h 2 & \h 3 & \h 4 & \h 5 \\
    \hline
     1 & & 1 \\
     2 & & 2 & 7 \\
     3 & & 3 & 18 & 15 \\
     4 & & 4 & 30 & 40 & 10 \\
     5 & & 5 & 40 & 60 & 20 & 1 \\
     6 & & 6 & 45 & 60 & 15 \\
     7 & & 7 & 42 & 35 \\
     8 & & 8 & 28 \\
     9 & & 9 \\
    10 & 10 %
\end{tabular}%
    }
\end{center}%
\end{minipage}%
\vspace{0.5\baselineskip}

Let us see how the matrices $c'{}^{i}_{jk}$ can be obtained from the 
matrices $c^{i}_{jk}$. Formula \eqref{e-37} can be written:
\[
    c'{}^{i}_{jk} = \sum_{f} c^{i-j}_{fk} \binom{j}{f}, 
    \qquad \mbox{i.e.,} \qquad
    \mathcal{C}'_{k} = \mathcal{C}_{k} \widetilde{\mathcal{P}} \tag{1,A}
\]
where $(\mathcal{C}'_{k})_{ij} = c'{}^{i+j}_{jk}$, 
$(\mathcal{C}_{~})_{~~} = c^{~}_{~~}$, and $(\mathcal{P})_{rs} = \binom{r}{s}$ 
is the Pascal matrix, so ($r$ rows, $s$ columns):
\[
(\widetilde{\mathcal{P}})_{rs} = \binom{s}{r} = 
\begin{tabular}{ c | r r r r r r r r }
            & \h 1 & \h 2 & \h 3 & \h 4 & \h 5 & \h 6 & \h 7 & $\cdots$ \\
    \hline
          1 & 1 & 2 & 3 & 4 &  5 &  6 &  7 & $\cdots$ \\
          2 &   & 1 & 3 & 6 & 10 & 15 & 21 & $\cdots$ \\
          3 &   &   & 1 & 4 & 10 & 20 & 35 & $\cdots$ \\
          4 &   &   &   & 1 &  5 & 16 & 35 & $\cdots$ \\
          5 &   &   &   &   &  1 &  6 & 21 & $\cdots$ \\
     \vdots &   &   &   &   &    &    &    & %
\end{tabular}
\]
We have succesively ($i$ rows, $j$ columns):
\vspace{0.5\baselineskip}
\begin{center}%
{$c'_{j}{}^{i}_{0} = $}
    {\footnotesize%
\begin{tabular}{ r | r r r r r r r r r r r }
       & \K 1 & \K 2 & \K 3 & \K 4 & \K 5 & \K 6 & \K 7 & \K 8 & \K 9 & \h 10 & \h 11 \\
    \hline
     1 & \\
     2 & 1 \\
     3 & & 2 \\
     4 & & & 3 \\
     5 & & & & 4 \\
     6 & & & & & 5 \\
     7 & & & & & & 6 \\
     8 & & & & & & & 7 \\
     9 & & & & & & & & 8 \\
    10 & & & & & & & & & 9 \\
    11 & & & & & & & & & & 10 \\
    12 & & & & & & & & & & & 11 %
\end{tabular}%
    }
\end{center}%
\vspace{0.5\baselineskip}
\begin{center}%
{$c'_{j}{}^{i}_{1} = $}
    {\footnotesize%
\begin{tabular}{ r | r r r r r r r r r r r }
       & \K 1 & \K 2 & \K 3 & \K 4 & \K 5 & \K 6 & \K 7 & \K 8 & \K 9 & \h 10 & \h 11 \\
    \hline
     1 & \\
     2 & \\
     3 & 1 \\
     4 & 1 & 2 \\
     5 & 1 & 4 & 3 \\
     6 & 1 & 4 & 9 & 4 \\
     7 & 1 & 4 & 12 & 16 & 5 \\
     8 & 1 & 4 & 12 & 28 & 25 & 6  \\
     9 & 1 & 4 & 12 & 32 & 55 & 36  & 7  \\
    10 & 1 & 4 & 12 & 32 & 75 & 96  & 49  & 8  \\
    11 & 1 & 4 & 12 & 32 & 80 & 156 & 154 & 64  & 9 \\
    12 & 1 & 4 & 12 & 32 & 80 & 186 & 294 & 232 & 81 & 10 %
\end{tabular}%
    }
\end{center}%
\vspace{0.5\baselineskip}
\begin{center}%
{$c'_{j}{}^{i}_{2} = $}
    {\footnotesize%
\begin{tabular}{ r | r r r r r r r r r r r }
       & \K 1 & \K 2 & \K 3 & \K 4 & \K 5 & \K 6 & \K 7 & \K 8 & \K 9 & \h 10 & \h 11 \\
    \hline
     1 & \\
     2 & \\
     3 & \\
     4 & \\
     5 & \\
     6 & & 1 \\
     7 & & 2 & 3 \\
     8 & & 3 & 9 & 6 \\
     9 & & 4 & 15 & 24 & 10  \\
    10 & & 5 & 21 & 48 & 50  & 15  \\
    11 & & 6 & 27 & 72 & 120 & 90  & 21  \\
    12 & & 7 & 33 & 96 & 200 & 255 & 147 & 28 %
\end{tabular}%
    }
\end{center}%
\vspace{0.5\baselineskip}
\begin{center}%
{$c'_{j}{}^{i}_{3} = $}
    {\footnotesize%
\begin{tabular}{ r | r r r r r r r r r r r }
       & \K 1 & \K 2 & \K 3 & \K 4 & \K 5 & \K 6 & \K 7 & \K 8 & \K 9 & \h 10 & \h 11 \\
    \hline
     1 & \\
     2 & \\
     3 & \\
     4 & \\
     5 & \\
     6 & \\
     7 & \\
     8 & \\
     9 & & & 1 \\
    10 & & & 3 & 4 \\
    11 & & & 6 & 16 & 10 \\
    12 & & & 10 & 36 & 50 & 20 %
\end{tabular}%
    }
\end{center}%
\vspace{0.5\baselineskip}

Finally, formula \eqref{e-43} can be written
\[
    c'_{k}{}^{i}_{j} = \sum_{f} \binom{j-f}{i-(k+2f)-(j-k)}%
                           \binom{j}{j-f} \binom{k-1}{f-1} ,
\]
which gives successively ($i$ rows, $j$ columns):

\vspace{0.5\baselineskip}
\noindent
\begin{minipage}{3in}%
\begin{center}%
{$c'_{0}{}^{i}_{j} = $}\\
\vspace{0.5\baselineskip}
    {\footnotesize%
\begin{tabular}{ r | r r r r r r r r r }
       & \h 1 & \h 2 & \h 3 & \h 4 & \h 5 & \h 6 & \h 7 & \h 8 & \h 9 \\
    \hline
     1 & 1 \\
     2 & 1 & 1 \\
     3 & & 2 & 1 \\
     4 & & 1 & 3 & 1 \\
     5 & & & 3 & 4 & 1 \\
     6 & & & 1 & 6 & 5 & 1 \\
     7 & & & & 4 & 10 & 6 & 1 \\
     8 & & & & 1 & 10 & 15 & 7 & 1 \\
     9 & & & & & 5 & 20 & 21 & 8 & 1 \\
    10 & & & & & 1 & 15 & 35 & 28 & 9 %
\end{tabular}%
    }
\end{center}%
\end{minipage}%
\hfill
\begin{minipage}{3in}%
\begin{center}%
{$c'_{1}{}^{i}_{j} = $}\\
\vspace{0.5\baselineskip}
    {\footnotesize%
\begin{tabular}{ r | r r r r r r r r r }
       & \h 1 & \h 2 & \h 3 & \h 4 & \h 5 & \h 6 & \h 7 & \h 8 & \h 9 \\
    \hline
     1 & \\
     2 & \\
     3 & 1 \\
     4 & & 2 \\
     5 & & 2 & 3 \\
     6 & & & 6 & 4 \\
     7 & & & 3 & 12 & 5 \\
     8 & & & & 12 & 20 & 6 \\
     9 & & & & 4 & 30 & 30 & 7 \\
    10 & & & & & 20 & 60 & 42 & 8 \\
\end{tabular}%
    }
\end{center}%
\end{minipage}%

\vspace{0.5\baselineskip}
\noindent
\begin{minipage}{3in}%
\begin{center}%
{$c'_{2}{}^{i}_{j} = $}\\
\vspace{0.5\baselineskip}
    {\footnotesize%
\begin{tabular}{ r | r r r r r r r r r }
       & \h 1 & \h 2 & \h 3 & \h 4 & \h 5 & \h 6 & \h 7 & \h 8 & \h 9 \\
    \hline
     1 & \\
     2 & \\
     3 & \\
     4 & 1 \\
     5 & & 2 \\
     6 & & 3 & 3 \\
     7 & & & 9 & 4 \\
     8 & & & 6 & 18 & 5 \\
     9 & & & & 30 & 50 & 6 \\
    10 & & & & 10 & 60 & 45 & 7 %
\end{tabular}%
    }
\end{center}%
\end{minipage}%
\hfill%
\begin{minipage}{3in}%
\begin{center}%
{$c'_{3}{}^{i}_{j} = $}\\
\vspace{0.5\baselineskip}
    {\footnotesize%
\begin{tabular}{ r | r r r r r r r r r }
       & \h 1 & \h 2 & \h 3 & \h 4 & \h 5 & \h 6 & \h 7 & \h 8 & \h 9 \\
    \hline
     1 &   \\
     2 &   \\
     3 &   \\
     4 &   \\
     5 & 1 \\
     6 &   & 2 \\
     7 &   & 4 &  3 \\
     8 &   &   & 12 &  4 \\
     9 &   &   & 20 & 24 &  5 \\
    10 &   &   &    & 40 & 40 & 6 %
\end{tabular}%
    }
\end{center}%
\end{minipage}%
\section*{Acknowledgements}

This work was partly done while in New York University on sabbatical leave. I 
thank Prof.\ Martin Pope of the Radiation and Solid State Laboratory and 
Prof.\ Peter Lax of the Courant Institute for their hospitality, and 
Prof.\ Frank Anger of the Department of Mathematics of the University of Puerto 
Rico at R\'\i o Piedras, for useful comments on the manuscript.


\begin{thebibliography}{99}

   \bibitem{Ad1910}
      Andr\'e, D. %
      \emph{Notice sur les travaux scientifiques}. %
      Gautier Villars, Paris, 1910.
      
   \bibitem{Bc1971}
      Berge, C. %
      \emph{Principles of Combinatorics}. %
      Math.\ in Sc.\ and Eng.\ \textbf{72}, %
      Academic Press, New York, 1971.
      
   \bibitem{Cl1974}
      Carlitz, L. %
      \emph{Fibonacci Notes, 1. Zero-one sequences and 
      Fibonacci numbers of higher order}. %
      Fibonacci Quart.\ \textbf{12} (1974), 1--10.

   \bibitem{Cl1976}
      \bysame\ %
      \emph{Fibonacci Notes, 5. Zero-one sequences again}. %
      Fibonacci Quart.\ (1977), 49--56.

   \bibitem{Cl1970}
      Comptet, L. %
      \emph{Analyse Combinatoire, vol.~I}. %
      Presses Universitaires de France, Paris, 1970, %
      pp.~32--34 and 93--96.
      
   \bibitem{Ki1953}
      Kaplansky, I. %
      \emph{Solution of the ``probl\`eme des m\'enages''}. %
      Bull.\ Amer.\ Math.\ Soc.\ \textbf{49} (1953), 784--785.

   \bibitem{Kde1968}
      Knuth, D. E. %
      \emph{The Art of Computer Programming, vol.~I: %
      Fundamental Algorithms}. %
      Addison-Wesley, Reading, MA, 1968.

   \bibitem{Mew1964}
      Montroll, E. W. %
      \emph{Lattice Statistics}, %
      in Applied Combinatorial Mathematics. %
      University of California Engineering and Physical Sciences %
      Extension Series. E.~F.~Beckenbach editor. %
      John Wiley, New York, 1964.

   \bibitem{Rj1958}
      Riordan, J. %
      \emph{An Introduction to Combinatorial Analysis}. %
      John Wiley, New York, 1958.

   \bibitem{Tcj1972}
      Thompson, C. J. %
      \emph{Mathematical Statistical Mechanics}. %
      Macmillan, New York, 1972.
      
\end{thebibliography}
\end{document}